\documentclass[12pt]{amsart}

\usepackage{amsfonts, amsmath, amsthm, amssymb} 
\usepackage{graphicx}
\usepackage{epstopdf}
\usepackage{psfrag}
\usepackage{pdfsync}

\newtheorem{pro}{Proposition}[section]
\newtheorem{thm}[pro]{Theorem}
\newtheorem{lem}[pro]{Lemma}
\newtheorem{clm}[pro]{Claim}

\newtheorem{cnj}[pro]{Conjecture}

\theoremstyle{definition}
\newtheorem{dfn}[pro]{Definition}
\newtheorem{ex}[pro]{Example}

\theoremstyle{remark}

\newtheorem*{note}{Note}

\newcommand{\thick}[1]{{\rm Thick}(#1)}
\newcommand{\thin}[1]{{\rm Thin}(#1)}

\newcommand{\VV}{\mathbb V}
\newcommand{\WW}{\mathbb W}
\newcommand{\CV}{\mathcal V}
\newcommand{\CW}{\mathcal W}

\newcommand{\Wup}{\CW _{\uparrow}}
\newcommand{\Wdown}{\CW _{\downarrow}}

\newcommand{\bdyup}{\partial _{\uparrow}}
\newcommand{\bdydown}{\partial _{\downarrow}}

\newcommand{\amlg}[1]{\mathcal A(#1)}

\title{Connected sums of unstabilized Heegaard splittings are unstabilized} 
\date{\today}

\author{David Bachman}
\address{Pitzer College\\1050 N. Mills Ave\\Claremont, CA 91711}
\email{bachman@pitzer.edu}

\begin{document}

\begin{abstract}
Let $M_1$ and $M_2$ be closed, orientable 3-manifolds. Let $H_i$ denote a Heegaard surface in $M_i$. We prove that if $H_1 \# H_2$ comes from stabilizing a lower genus splitting of $M_1 \# M_2$ then either $H_1$ or $H_2$ comes from stabilizing a lower genus splitting. If $H_i$ and $G_i$ are non-isotopic Heegaard surfaces in $M_i$, and $H_i$ is unstabilized, then we show $H_1 \# H_2$ is not isotopic to $G_1 \# G_2$ in $M_1 \# M_2$. The former result answers a question of C. Gordon (Problem 3.91 from \cite{kirby:97}).
\end{abstract}
\maketitle

\pagestyle{myheadings}
\markboth{DAVID BACHMAN}{CONNECTED SUMS OF UNSTABILIZED HEEGAARD SPLITTINGS}

\noindent
Keywords: Heegaard Splitting, Incompressible Surface

\section{Introduction}

Suppose $M_1$ and $M_2$ are closed, orientable 3-manifolds and $H_i$ is a Heegaard surface in $M_i$. Then one can form the connected sum $H=H_1 \# H_2$ in the 3-manifold $M=M_1 \#M_2$ to obtain a new Heegaard surface. If $H_1$, say, came from stabilizing some lower genus Heegaard surface in $M_1$ then it immediately follows that $H$ comes from stabilizing a lower genus Heegaard surface in $M$. In 1997 C. McA. Gordon conjectured that the converse must also be true (see Problem 3.91 from \cite{kirby:97}):

\begin{cnj}[Gordon's Conjecture]
Suppose $H_i$ is a Heegaard surface in $M_i$, for $i=1,2$. If $H_1 \# H_2$ is a stabilized Heegaard surface in $M_1 \# M_2$ then either $H_1$ or $H_2$ is stabilized. 
\end{cnj}

This paper contains a proof of Gordon's conjecture (Theorem \ref{t:main}). In a previous version of this paper a proof of the conjecture was announced with the additional assumption that $M_1$ and $M_2$ are irreducible. At the same time Ruifeng Qui announced a complete proof of Gordon's conjecture \cite{qiu:04}. The assumption of irreducibility is not used in the present version of this paper. In addition, we prove the following:

\medskip
\noindent {\bf Theorem \ref{t:IsotopyTheorem}.} {\it Let $M_1$ and $M_2$ be two closed, orientable 3-manifolds. Suppose $H_i $ and $G_i$ are non-isotopic Heegaard splittings of $M_i$ and $H_i$ is unstabilized. Then $H_1 \# H_2$ is not isotopic to $G_1 \# G_2$ in $M_1 \# M_2$.}
\medskip

The surfaces $H_1$ and $H_2$ of Gordon's conjecture together form a structure called a {\it generalized Heegaard splitting (GHS)}. Loosely speaking, a GHS $H$ is a pair of sets of surfaces, $\thick{H}$ and $\thin{H}$, such that each element of $\thick{H}$  is a Heegaard surface for some component of $M-\thin{H}$ (see Section \ref{s:GHSdefinition} for a more precise definition). 

We can transform one GHS into another by a process called {\it weak reduction}. This can be done whenever there are disjoint compressing disks on opposite sides of some thick surface  (see Section \ref{s:WRdefinition}). So, for example, if some thick surface of a GHS comes from stabilization or connected sum then there is a weak reduction for the entire GHS. Definition \ref{d:GHScomplexity} gives a simple complexity for GHSs under which weak reduction represents a decrease. 

Once we can relate GHSs by weak reduction we can start examining an entire {\it Sequence of GHSs (SOG)}. This is defined to be a sequence $\{H^i\}$ such that for each $i$ either the GHS $H^{i+1}$ or the GHS $H^i$ can be obtained from the other by weak reduction. Now, given fixed GHSs $H$ and $H'$ and an SOG which connects them, one can ask if there is some sense in which there is a ``more efficient" SOG which connects them. In Section \ref{s:SOGdefinition} we define several ways to find such a new SOG. Any SOG obtained by one of these operations is said to have been obtained by a {\it reduction}. If a simpler SOG cannot be found then the given one is said to be {\it irreducible}. 

Section \ref{s:SOGdefinition}  concludes with a crucial result about irreducible SOGs. This is given by Lemma \ref{l:maximalGHS}, which states that the thick surfaces of the maximal GHSs of an irreducible SOG satisfy one of two combinatorial conditions. These conditions are called {\it strong irreducibility} and {\it criticality}. Strongly irreducible Heegaard splittings were introduced by Casson and Gordon in \cite{cg:87}.  Criticality was introduced by the author in \cite{crit},  although the definition given here in Section \ref{s:SICrit} is considerably simpler. By Lemma \ref{l:SIorCritImpliesIrreducible} any GHS whose thick surfaces are strongly irreducible or critical must be of an irreducible 3-manifold. This lemma relies on a deep result about the intersections of strongly irreducible and critical surfaces with incompressible surfaces (Lemma \ref{l:EssentialIntersection}).

Our proof of Gordon's conjecture (Theorem \ref{t:main}) begins with the construction of a SOG $\{H^i\}_{i=1}^{n}$ as follows. Each element of the set $\thick{H^1}$ is a Heegaard surface of an irreducible 3-manifold. For some $k \ge 1$ the set $\thick{H^k}$ has two elements, namely $H_1$ and $H_2$, where $H_i$ is the connected sum of a subset of $\thick{H^1}$ for $i=1,2$. The set $\thick{H^{k+1}}$ consists of a single element, $H=H_1 \# H_2$. We assume that $H$ destabilizes to a Heegaard surface $G$, and define $H^{k+2}$ to be the GHS with $\thick{H^{k+2}}=\{G\}$. Finally, $H^n$ is a GHS where again each element of $\thick{H^n}$ is a Heegaard surface in an irreducible 3-manifold, and $G$ is the connected sum of these surfaces. It follows that the sum of the genera of the thick surfaces of $H^n$ is strictly less than the sum of the genera of the thick surfaces of $H^1$. 

The SOG $\{H^i\}$ thus defined has a single maximal GHS, $H^{k+1}$. Since this is not a GHS of an irreducible 3-manifold, there must be a more efficient SOG from $H^1$ to $H^n$. In other words, there is some reduction for $\{H^i\}$. Reducing as much as possible yields a SOG in which every element is a GHS of an irreducible 3-manifold. Finally, it follows from the assumption that $H_1$ and $H_2$ are unstabilized that there are no destabilizations in this final SOG. This gives our contradiction, as it implies that the sum of the genera of the thick sufaces of $H^1$ is equal to the sum of the genera of the thick surfaces of $H^n$.

The author would like to thank Saul Schleimer for helpful conversations regarding the proof of Claim \ref{c:PhiExists}, and the referee for several extremely helpful suggestions, including many of the examples given in the paper. 

\section{Definitions.}

\subsection{Essential loops, disks and spheres}
A 2-sphere in a 3-manifold which does not bound a 3-ball is called {\it essential}. If a manifold does not contain an essential 2-sphere then it is referred to as {\it irreducible}. 

A loop on a surface is called {\it essential} if it does not bound a disk in the surface.  An arc which is properly embedded in a surface $F$ is {\it essential} if it does not cobound, with a subarc of $\partial F$, a subdisk of $F$. 

Suppose $F$ is a surface embedded in a 3-manifold $M$, $D$ is a disk in $M$, and $D \cap F=\partial D$. There is an embedding $h:D \times I \to M$ such that $h(D \times \{\frac{1}{2}\})=D$ and $h(D \times I) \cap F=h(\partial D \times I)$. To {\it surger $F$ along $D$} is to remove $h(\partial D \times I)$ from $F$ and replace it with $h(D \times \partial I)$. We denote the result of such a surgery as $F/D$. As a shorthand we will denote $(F/D)/E$ as $F/DE$ when $D$ and $E$ are disjoint disks with boundary on $F$.

If $\partial D$ is an essential loop on $F$ then $D$ is referred to as a {\it compressing disk for $F$}, and surgery along $D$ is referred to as {\it compression}. A surface $F$ is said to be {\it incompressible} if there are no compressing disks for $F$. A properly embedded disk in a 3-manifold $M$ is {\it essential} if it is a compressing disk for $\partial M$. 

Now suppose $F$ is a properly embedded surface in a 3-manifold $M$ with boundary, and $D$ is a disk such that $\partial D=\alpha \cup \beta$, $F \cap D=\alpha$ is an arc on $F$, and $D \cap \partial M=\beta$. Then there is an embedding $h:D \times I \to M$ such that $h(D \times \{\frac{1}{2}\})=D$ and $h(D \times I) \cap F=h(\alpha \times I)$ and $h(D \times I) \cap \partial M=h(\beta \times I)$. To {\it surger $F$ along $D$} is to remove $h(\alpha \times I)$ from $F$ and replace it with $h(D \times \partial I)$. If, furthermore, $\alpha$ is an essential arc on $F$ then $D$ is referred to as a {\it $\partial$-compressing disk for $F$}, and surgery along $D$ is referred to as {\it $\partial$-compression}. A surface $F$ is said to be {\it $\partial$-incompressible} if there are no $\partial$-compressing disks for $F$. 

If $M_1$ and $M_2$ are $n$-manifolds then the connected sum, denoted $M_1 \# M_2$, is constructed as follows. First, obtain $M_i^*$ by removing the interior of an $n$-ball from the interior of $M_i$. Each $M_i^*$ will thus have a new $(n-1)$-sphere boundary component, $S_i$. The manifold $M_1 \# M_2$ is then obtained by identifying $S_1$ with $S_2$. The image of $S_i$ in $M_1 \# M_2$ is referred to as the {\it summing sphere}. Note that the summing sphere is essential if and only if neither $M_i$ is an $n$-sphere. 

%I don't think I need this anymore:
%\begin{lem}
%\label{l:UniqueEssential}
%Suppose $M=M_1 \# M_2$, where each $M_i$ is an irreducible 3-manifold. Let $S$ denote the summing sphere. If $S'$ is any essential sphere in $M$ then $S'$ is isotopic to $S$. 
%\end{lem}

%\begin{proof}
%As above $M=M_1^* \cup M_2^*$, where $M_i^*$ is obtained from $M_i$ by removing the interior of a 3-ball.  Isotope $S'$ so that $|S \cap S'|$ is minimal. We claim that now $S \cap S'=\emptyset$. Let $D'$ denote an innermost subdisk of $S'$ bounded by a loop of $S \cap S'$. Then $D'$ lies in either $M_1^*$ or $M_2^*$. The loop $\partial D'$ divides the sphere $S$ into two disks. The irreducibility of $M_i$ then implies that one of these disks, together with $D'$, forms a sphere which bounds a ball in $M_i^*$. Call this disk $D$. We now use the ball bounded by $D \cup D'$ in $M_i^*$ to guide an isotopy of $S'$, removing at least one loop of $S \cap S'$. 

%If $S\cap S'=\emptyset$ then the result follows as $S'$ lies in a submanifold of $M$ homeomorphic to a punctured, irreducible 3-manifold with boundary $S$. 
%\end{proof}

\subsection{Heegaard splittings}
\label{s:Heegaard}

\begin{dfn}
A {\it compression body} is a 3-manifold which can be obtained by starting with some closed, orientable, connected surface, $H$, forming the product $H \times I$, attaching some number of 2-handles to $H \times \{1\}$, and capping off  all resulting 2-sphere boundary components with 3-balls. The boundary component $H \times \{0\}$ is referred to as $\partial _+$. The rest of the boundary is referred to as $\partial _-$.  
\end{dfn}

\begin{dfn}
A {\it Heegaard splitting} of a 3-manifold $M$ is an expression of $M$ as a union $\mathcal V \cup _H \mathcal W$, where $\mathcal V$ and $\mathcal W$ are compression bodies that intersect in an oriented surface $H=\partial _+ \mathcal V=\partial _+ \mathcal W$. If $\mathcal V \cup _H \mathcal W$ is a Heegaard splitting of $M$ then we say $H$ is a {\it Heegaard surface}.
\end{dfn}

\begin{note}
The assumption that $H$ is oriented in the above definition will play an important role. For example, $L(p,q)$, where $q \ne  \pm 1\ {\rm mod}\ p$, contains a pair of non-isotopic Heegaard tori \cite{bo:83}, but as unoriented surfaces they are isotopic. In contrast, $S^3$ contains a unique Heegaard torus \cite{waldhausen:68}, namely the boundary of a regular neighborhood of an unknotted loop. Hence, we may talk about the {\it standard} genus one Heegaard surface in $S^3$.
\end{note}

\begin{dfn}
\label{d:ConnectedSum}
Suppose $M_1$ and $M_2$ are 3-manifolds and $H_i$ is a Heegaard surface in $M_i$. Recall that the first step in defining the connected sum $M_1 \# M_2$ is to remove the interior of a ball $B_i$ from $M_i$, resulting in a new 2-sphere boundary component $S_i$ of the punctured 3-manifold $M_i^*$. If $B_i$ is chosen to meet $H_i$ in a disk then $H_i^*=H_i \cap M_i^*$ will divide $S_i$ into disks $D_i$ and $D_i'$. Now the identification of $S_1$ with $S_2$ can be done in two ways; $D_1$ is glued to $D_2$ or to $D_2'$. However, only one such identification will make the orientation of $H_1^*$ agree with that of $H_2^*$. When this identification is used the surface $H_1^* \cup H_2^*$ is referred to as the connected sum $H_1 \# H_2$ of $H_1$ and $H_2$. 
\end{dfn}

\begin{ex}
\label{e:LensSpaces}
Let $H_1$ and $H_2$ be Heegaard tori in $L(p,q)$. Let $\overline{H_1}$ denote $H_1$ with the opposite orientation. Then $H_1 \# H_2$ and $\overline{H_1} \# H_2$ are non-isotopic genus two Heegaard surfaces in $L(p,q) \# L(p,q)$, even as unoriented surfaces \cite{engmann:70} (see also \cite{birman:74}).
\end{ex}

\begin{dfn}
\label{d:stabilization}
A {\it stabilization} of a Heegaard surface $H$ is a new Heegaard surface which is the connected sum of $H$ with the standard genus one Heegaard surface in $S^3$. 
\end{dfn}

%A stabilization $G$ of a Heegaard surface $H$ is completely characterized by the presence of a pair of compressing disks $(D,E)$ on opposite sides of $G$ which meet in a point, such that compressing $G$ along either $D$ or $E$ yields a surface isotopic to $H$.

%\begin{ex}
%Let $H$ be a Heegaard surface in a punctured lens space, $L$. Then $H$ and the Heegaard surface $\overline{H}$ with opposite orientation are not isotopic. However, if we stabilize $H$ once then we obtain a Heegaard surface that is isotopic to one stabilization of $\overline{H}$. We demonstrate this as follows.

%Construct $L$ by starting with the product $T^2 \times I$ and gluing a solid torus to $T^2 \times \{0\}$ and a punctured solid torus to $T^2 \times\{1\}$. Let $D$ be a disk  in $T^2$ and let $T^*=T^2-int(D)$. The surfaces $T^2 \times \{0\}$ and $T^2 \times \{1\}$ are Heegaard tori in $L$. Orient these surfaces oppositely. Now note that $(T^* \times \{0\}) \cup (\partial D \times I) \cup (T^* \times \{1\})$ is a common stabilization.
%\end{ex}

%If the Heegaard surface bounds compression bodies $\CV$ and $\CW$ then another way to define a stabilization is by ``tunneling" a 1-handle out of $W$ and attaching it to $\partial _+ \CW$. If this is done in such a way so as to make the definition symmetric in $\CV$ and $\CW$ then one arrives at a stabilization. 
%The Reidemeister-Singer theorem states that given any two Heegaard surfaces, $H$ and $H'$, there is always a stabilization of $H$ which is isotopic to a stabilization of $H'$ (see, for example, \cite{am:90}). 

\section{Strong Irreducibility and Criticality}
\label{s:SICrit}

%Can simplify this section slightly by assuming that all surfaces are connected?

The main technical tools of this paper are {\it strongly irreducible} \cite{cg:87} and {\it critical} \cite{crit} surfaces. Both strong irreducibilty and criticality are combinatorial conditions satisfied by the compressing disks for a Heegaard surface.
%thick level of a GHS for which there are no weak reductions or destabilizations. Such GHSs correspond to handle structures in which the 2-handles are always added before the 1-handles, when possible (see \cite{st:94}). The definition of criticality is also a combinatorial condition satisfied by the compressing disks for thick levels of some GHSs. These GHSs correspond to maxima of SOGs in which amalgamations and stabilizations come before weak reductions and destabilizations, when possible. 

\begin{dfn}
\label{d:reducingpair}
Let $\mathcal V \cup _H \mathcal W$ be a Heegaard splitting of a 3-manifold $M$. Then we say the pair $(V,W)$ is a {\it reducing pair} for $H$ if $V$ and $W$ are disjoint compressing disks on opposite sides of $H$.
\end{dfn}

\begin{dfn}
A Heegaard surface is {\it strongly irreducible} if it is compressible to both sides but has no reducing pairs.
\end{dfn}

%Under new dfn the genus 1 splitting of S^3 is SI. Do I care?

%Maybe take out assumption that it is compressible, which would imply the SI->Inc or compressible to both sides.

\begin{dfn}
\label{d:critical}
Let $H$ be a Heegaard surface in some 3-manifold which is compressible to both sides. The surface $H$ is {\it critical} if the set of all compressing disks for $H$ can be partitioned into subsets $C_0$ and $C_1$ such that
\begin{enumerate}
	\item For each $i=0,1$ there is at least one reducing pair $(V_i,W_i)$, where $V_i, W_i \in C_i$.
	\item If $V \in C_0$ and $W \in C_1$ then $(V,W)$ is not a reducing pair. 
\end{enumerate}
\end{dfn}

Definition \ref{d:critical} is significantly simpler, and slightly weaker, than the one given in \cite{crit}.  In other words, anything that was considered critical in \cite{crit} is considered critical here as well. 
%Is this really still true?
Hence, a result such as Theorem 7.1 of \cite{crit} still holds. This result says that in a non-Haken 3-manifold the minimal genus common stabilization of any pair of non-isotopic, unstabilized Heegaard splittings is critical. The basic idea of the proof is as follows: suppose $H_0$ and $H_1$ are non-isotopic Heegaard splittings in a 3-manifold $M$ which are isotopic to a surface $K$ after one stabilization. Another way to say this is that there are reducing pairs $(V_0,W_0)$ and $(V_1,W_1)$ representing destabilizations
%Ack! Reducing pairs no longer represent destabilizations.
 of $K$ that lead to $H_0$ and $H_1$. We then show that either we can use the disks $V_i$ and $W_i$ to create a partition of the compressing disks for $K$ that satisfies the conditions of Definition \ref{d:critical}, or there is an incompressible surface in $M$.

\begin{ex}
\label{e:smallSFS}
Let $M$ be a Seifert fibered space which fibers over the sphere with three exceptional fibers.  There are three vertical splittings $H$, $H'$, and $H''$ of $M$ and these are generally not isotopic (see \cite{ms:98} for the relevant definitions).  Let $K$ be the genus three splitting which is the common stabilization of these three.  Since $M$ is non-Haken, it follows from Theorem 7.1 of \cite{crit} that $K$ is critical. 
\end{ex}

Example \ref{e:smallSFS} shows that the partition of disks into just two sets in Definition \ref{d:critical} is a bit misleading. For a critical Heegaard surface in a non-Haken 3-manifold one can make a partition with a set for each distinct destabilization. This point is made more explicit in the definition of criticality given in \cite{crit}. Somehow the useful results (such as Lemma \ref{l:EssentialIntersection} below) about critical surfaces only require a partition with at least two sets.  This is one reason for the more streamlined definition given here.

\begin{lem}
\label{l:EssentialIntersection}
Suppose $H$ is an incompressible,  strongly irreducible, or critical Heegaard surface in a 3-manifold $M$ and $S$ is an essential sphere or disk in $M$. Then there is an essential surface $S'$, obtained from $S$ by surgery, such that $S' \cap H=\emptyset$.
\end{lem}

\begin{proof}
We break the proof into three cases, depending on whether $H$ is incompressible, strongly irreducible, or critical.

\medskip

\noindent {\bf Case 1.} {\it $H$ is incompressible.} The incompressible case follows from a standard innermost disk argument. We leave the proof to the reader.

\medskip

\noindent {\bf Case 2.}  {\it $H$ is strongly irreducible.} The strongly irreducible case also follows from standard arguments. For completeness we present a proof here. 

Let $H$ be a strongly irreducible Heegaard surface, separating $M$ into $\CV$ and $\CW$. By the definition of strong irreducibility there are  compressing disks $V \subset \CV$ and $W \subset  \CW$ for $H$. A neighborhood of $V$ is homeomorphic to a 3-ball, and hence we can define an isotopy which pushes $S$ off $V$. Similarly, we can define an isotopy of $S$ which pushes $S$ off of $W$. Putting such isotopies together gives us an isotopy $S_t$ ({\it i.e.} a map $\gamma : S \times [-1,1] \to M$ where $S_t=\gamma(S,t))$ such that
\begin{itemize}
	\item $S_{-1} \cap V=\emptyset$
	\item $S_0=S$
	\item $S_{1} \cap W=\emptyset$
\end{itemize}
Let $t_0=-1$, $\{t_i\}_{i=1}^{n-1}$ be the set of points in $[-1,1]$ where $S_t$ is not transverse to $H$, and $t_n=1$. These points break $[-1,1]$ up into subintervals, which we now label. If there is a $t \in (t_i,t_{i+1})$ such that $H \cap S_t$ contains the boundary of a compressing disk for $H$ in $\CV$ then this interval gets the label $\VV$. Similarly, if there is a $t \in (t_i,t_{i+1})$ such that $H \cap S_t$ contains the boundary of a compressing disk for $H$ in $\CW$ then this interval gets the label $\WW$.

We now present several claims which produce the desired result.

\begin{clm}
The interval $(t_0,t_1)$ is either labeled $\VV$ or has no label. Similarly, the interval $(t_{n-1},t_n)$ is labeled $\WW$ or has no label.
\end{clm}

\begin{proof}
For $t$ near $-1$ the surface $S_t$ is disjoint from $V$. Suppose the interval containing $t$ is labeled $\WW$. Then there is a loop $\alpha \subset H \cap S_t$ which bounds a compressing disk $W'$ for $H$ in $\CW$. But then $\partial W' \cap \partial V=\emptyset$, contradicting the strong irreducibility of $H$. A symmetric argument completes the proof.
\end{proof} 

\begin{clm}
No interval has both labels.
\end{clm}

\begin{proof}
If, for some $t$, there are loops of $S_t \cap H$ that bound disks in $\CV$ and $\CW$ then we immediately contradict the strong irreducibility of $H$.
\end{proof}

\begin{clm}
\label{c:Gabai}
An interval labeled $\VV$ cannot be adjacent to an interval labeled $\WW$. 
\end{clm}

\begin{proof}
Suppose $t_i$ is the intersection point of adjacent intervals with different labels. Let $t_-=t_i-\epsilon$ and $t_+=t_i+\epsilon$. As $H$ is orientable the loops of $S_{t_-} \cap H$ can be made disjoint (on $H$) from the loops of $S_{t_+} \cap H$ (see, for example, Lemma 4.4 of \cite{gabai:87}). Hence, if $S_{t_-} \cap H$ contains the boundary of a compressing disk for $H$ in $\CV$ it can be made disjoint from the boundary of a compressing disk in $\CW$ that is contained in $S_{t_+} \cap H$. This again contradicts the strong irreducibility of $H$.
\end{proof}

Following these claims we conclude there is an unlabeled interval. Henceforth we assume $t$ is in such an interval. We now claim that every loop of $S_t \cap H$ is inessential on $H$. Suppose this is not the case. Let $\alpha$ be a loop of $S_t \cap H$ which is innermost on $S_t$ among all such loops that are also essential on $H$.

\begin{clm}
The loop $\alpha$ bounds a compressing disk for $H$. 
\end{clm}

\begin{proof}
As $S$ is a sphere or disk the loop $\alpha$ bounds a subdisk $A$ of $S_t$ so that all curves of $int(A) \cap H$ are inessential on $H$. If the interior of $A$ misses $H$ then $A$ is a compressing disk for $H$. If not then there is some loop $\beta$ where the interior of $A$ meets $H$. By assumption $\beta$ is inessential on $H$. Hence, $\beta$ bounds a subdisk $B$ of $H$. Let $\gamma$ denote a loop of $A \cap B$ which is innermost on $B$. Then $\gamma$ bounds a subdisk $C$ of $H$ whose interior is disjoint from $A$. We now use $C$ to surger $A$, removing one loop of $A \cap H$. Continuing in this way we may remove all loops of $A \cap H$, besides $\partial A$. The resulting disk is a compressing disk for $H$.
\end{proof}
 
By the previous claim the loop $\alpha$ bounds a compressing disk $D$ for $H$. This contradicts the fact that $t$ lies in an unlabeled interval. We conclude $S_t \cap H$ consists of loops that are inessential on both surfaces. A standard innermost disk argument now completes the proof of the strongly irreducible case of Lemma \ref{l:EssentialIntersection}.

\medskip

\noindent {\bf Case 3.} {\it $H$ is critical.} In this case Lemma \ref{l:EssentialIntersection} essentially follows from Theorem 5.1 of \cite{crit}. As we are using a slightly weaker definition of the term ``critical" here, we reproduce the proof. Henceforth we will assume that $H$ is a critical Heegaard surface which separates $M$ into compression bodies $\CV$ and $\CW$. Let $C_0$ and $C_1$ be the sets in Definition \ref{d:critical}. By definition there are compressing disks $V_i \subset \CV$ and $W_i \subset \CW$ where $(V_i,W_i)$ is a reducing pair and $V_i,W_i \in C_i$. 
%By Lemma \ref{l:CritImpliesWRpairs} we may assume that the disks $V_i$ and $W_i$ are disjoint. 

%If, for $i=0$ or $1$,  $V_i \cap W_i$ is a point then, if necessary, perform a small isotopy of $S$ to make $S \cap V_i \cap W_i=\emptyset$.

The proof is in several stages. First, we construct a map $\Phi$ from $S \times D^2$ into $M$. We then use $\Phi$ to break up $D^2$ into regions and label them in such a way so that if any region remains unlabelled then the conclusion of the lemma follows. Finally, we construct a map from $D^2$ to a labelled 2-complex $\Pi$ which has non-trivial first homology, and show that if there is no unlabelled region then the induced map on homology is nontrivial, a contradiction. This general strategy is similar to that used in \cite{rs:96}, although the details have little in common. 

\subsection{Constructing the map, $\Phi :S \times D^2 \rightarrow M$.}\

We begin by defining a two parameter family of surfaces in $M$ isotopic to $S$. For any map  $\Phi :S \times D^2 \rightarrow M$ we let $S_x$ denote the image of $\Phi(S,x)$. 

\begin{clm}
\label{c:PhiExists}
Let $\theta _0$ and $\theta _1$ be distinct points on $\partial D^2$. Let $U$ and $L$ be arcs of $\partial D^2$ such that $\partial D^2=U \cup L$, where $U \cap L=\theta_0 \cup \theta _1$. There is a continuous map $\Phi:S \times D^2 \rightarrow M$ such that 
\begin{itemize}
	\item for all $x \in D^2$ the surface $S_x$ is embedded,
	\item for $i=0,1$ the surface $S _{\theta_i} $ is disjoint from both $V_i$ and $W_i$,
	\item for each $\theta \in U$ the surface $S _{\theta}$ is disjoint from at least one compressing disk for $H$ in $\CV$, and 
	\item for each $\theta \in L$ the surface $S _{\theta}$ is disjoint from at least one compressing disk for $H$ in $\CW$.
\end{itemize} 
\end{clm}

\begin{proof}
We start by inductively defining a sequence of compressing disks for $H$, $\{V^i\}_{i=0}^{n}$, such that $V^i \cap V^{i+1} =0$, for all $i$ between 0 and $n-1$, $V^0=V_0$, and $V^n=V_1$. 

\begin{enumerate}
	\item Define $V^0=V_0$.
	\item Let $V^i$ denote the last disk defined, and suppose there is a simple closed curve in $V_1 \cap V^i$. Let $v$ denote an innermost subdisk of $V_1$ bounded by a loop of $V_1 \cap V^i$. Now surger $V^i$ along $v$. The result is a disk and a sphere. Throw away the sphere and denote the disk as $V^{i+1}$. Note that $|V_1 \cap V^{i+1}|<|V_1 \cap V^i|$.
	\item Let $V^i$ denote the last disk defined, and suppose there are only arcs in $V_1 \cap V^i$. Let $v$ denote an outermost subdisk of $V_1$ cut off by an arc of $V_1 \cap V^i$. Now surger $V^i$ along $v$. The result is two disks, at least one of which is a compressing disk for $H$. Call such a disk $V^{i+1}$. Again, note that $|V_1 \cap V^{i+1}|<|V_1 \cap V^i|$.
	\item If neither of the previous two cases apply then we have arrived at a disk $V^{n-1}$ such that $V^{n-1} \cap V_1=\emptyset$. Now let $V^n=V_1$ and we are done.
\end{enumerate}

We can apply a symmetric construction to produce a sequence of compressing disks, $\{W^j\}_{j=0}^{m}$, such that $W^0=W_0$, $W^m=W_1$, and $W^j \cap W^{j+1} =\emptyset$, for all $j$ between 0 and $m-1$.

        \begin{figure}[htbp]
        \psfrag{0}{$\theta_0$}
        \psfrag{1}{$\theta_1$}
        \psfrag{a}{$V^0$}
        \psfrag{b}{$V^2$}
        \psfrag{c}{$V^1$}
        \psfrag{x}{$W^0$}
        \psfrag{y}{$W^1$}
        \psfrag{U}{$U$}
        \psfrag{L}{$L$}
        \psfrag{H}{$H$}
        \psfrag{S}{$S$}
        \vspace{0 in}
        \begin{center}
        \includegraphics[width=5 in]{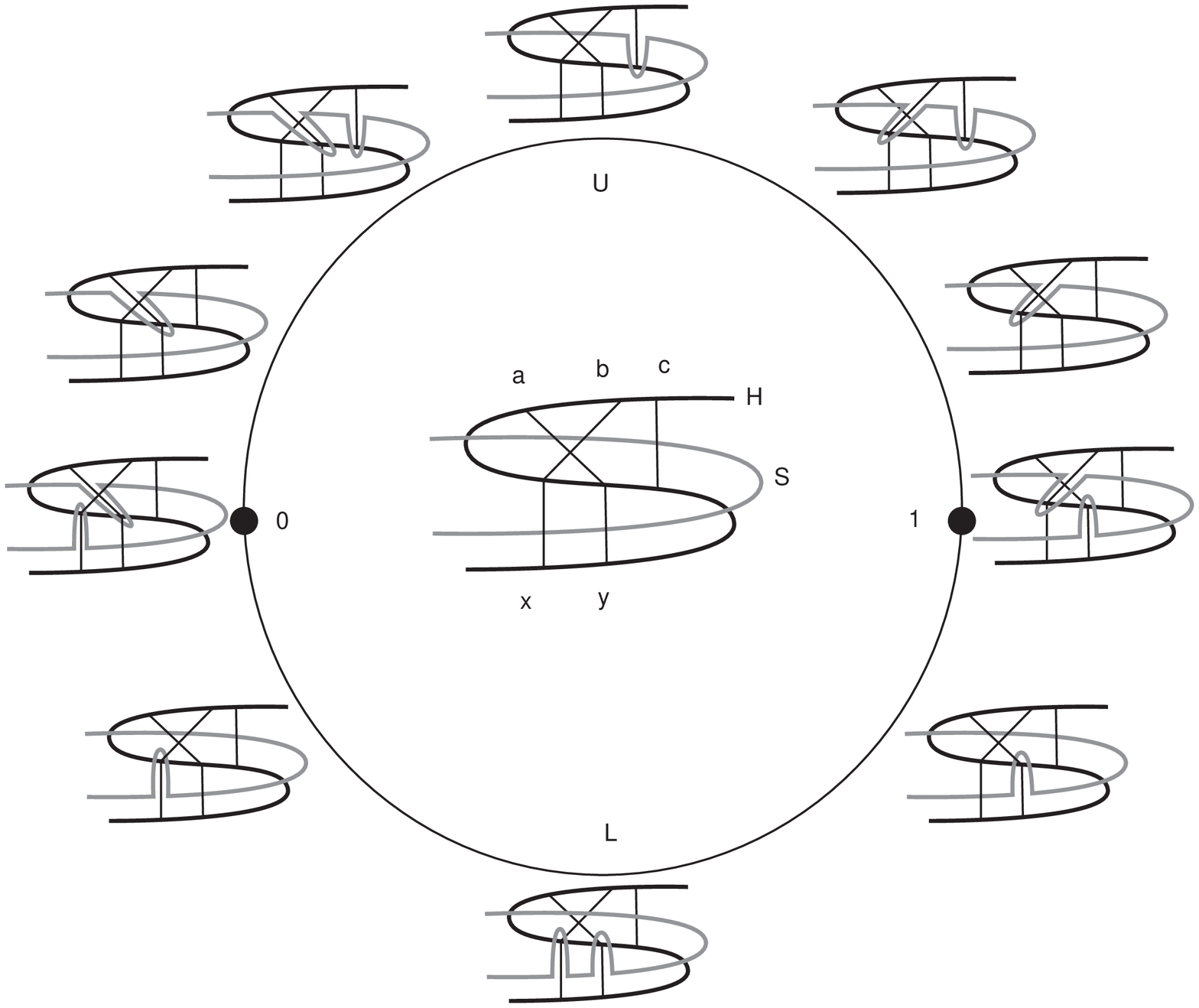}
        \caption{}
        \label{f:BigPhi}
        \end{center}
        \end{figure}

The map $\Phi$ can now be described as follows. See Figure \ref{f:BigPhi}. For $x$ at the center of $D^2$ the surface $S_x$ is identical to $S$. Near $\theta _0$ the surface $S_x$ is disjoint from both $V^0$ and $W^0$. As $x$ progresses along $U$ toward $\theta _1$ the surface $S_x$ ceases to be disjoint from $W^0$, but becomes disjoint from $V^1$. Progressing further the surface $S_x$ ceases to be disjoint from $V^0$, but becomes disjoint from $V^2$. This continues until $x$ gets to $\theta _1$, when $S_x$ is disjoint from both $V^n$ and $W^m$. This is illustrated for $n=2$ and $m=1$ in Figure \ref{f:BigPhi}.

To rigorously define $\Phi$ requires a considerable amount of further work (and, unfortunately, notation). For each $i$ between $0$ and $n$ let $A^i$ denote a neighborhood of the disk $V^i$. For each $j$ between $0$ and $m$ let $B^j$ denote a neighborhood of the disk $W^j$. Because $V^i \cap V^{i+1}=W^j \cap W^{j+1}=V^0 \cap W^0=V^n \cap W^m=\emptyset$ we may assume $A^i \cap A^{i+1}=B^j \cap B^{j+1}=A^0 \cap B^0=A^n \cap B^m=\emptyset$. 

For each $i$ between 0 and $n$, let $\gamma ^i :\overline{A^i} \times I \rightarrow \overline{A^i}$ be an isotopy which pushes $S$ off of $V^i$. In other words, 
%if $\gamma ^i _t(x)$ is short-hand for $\gamma ^i (x,t)$, then  
$\gamma ^i  (x,0)=x$ for all $x \in \overline{A^i}$, $\gamma ^i (S,1) \cap V^i =\emptyset$, and $\gamma ^i (x,t)=x$ for all $x \in \partial \overline{A^i}$. Similarly, for each $j$ between 0 and $m$, let $\delta ^j$ be an isotopy which pushes $S$ off of $W^j$, inside $B^j$. 

Choose $n$ pairs of circular arcs centered on points of $U$, and $m$ pairs of arcs centered on points of $L$, ``linked" as in Figure \ref{f:SaulPhi}\footnote{Figure suggested by Saul Schleimer.}. For the $i$th pair of arcs chosen, centered on a point of $U$, define $f^i:D^2 \to [0,1]$ to be the continuous function depicted in Figure \ref{f:SaulPhi2}. Let $g^j:D^2 \to [0,1]$ be the function similarly defined for the $j$th pair of arcs centered on a point of $L$.

        \begin{figure}[htbp]
        \psfrag{0}{$\theta_0$}
        \psfrag{1}{$\theta_1$}
        \vspace{0 in}
        \begin{center}
        \includegraphics[width=3 in]{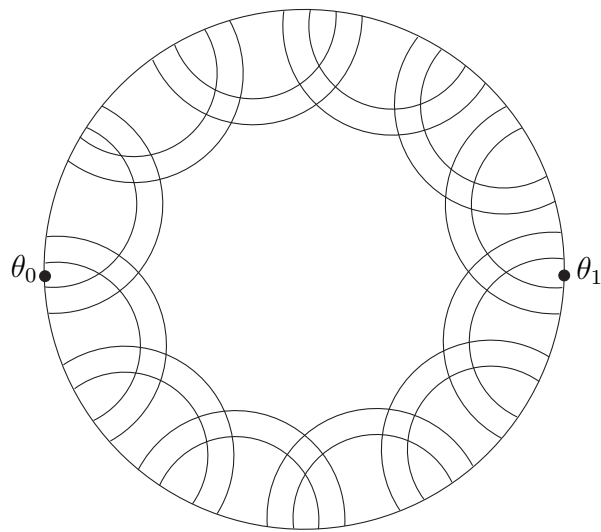}
        \caption{``Linked" pairs of circular arcs in $D^2$.}
        \label{f:SaulPhi}
        \end{center}
        \end{figure}

        \begin{figure}[htbp]
        \psfrag{0}{$\theta_0$}
        \psfrag{1}{$\theta_1$}
        \psfrag{f}{$f^i=1$}
        \psfrag{g}{$f^i=0$}
        \vspace{0 in}
        \begin{center}
        \includegraphics[width=3 in]{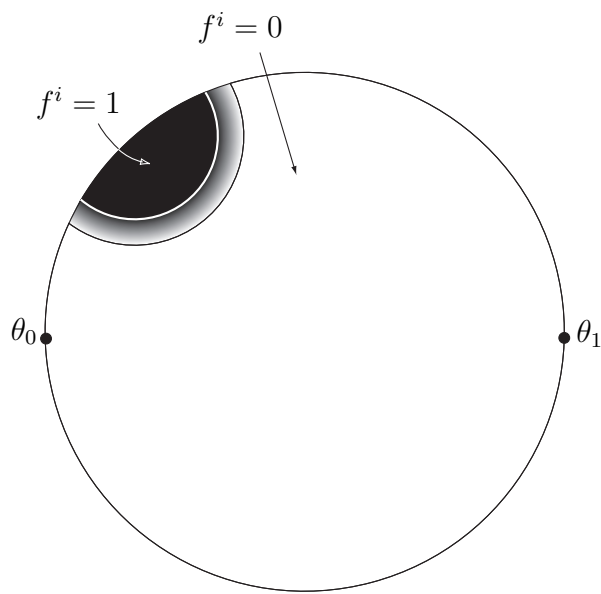}
        \caption{In the black region $f^i$ takes on the value 1. In the white region $f^i=0$. The shading between these two regions is meant to indicate that between the arcs the function $f^i$ continuously varies from $0$ to $1$.}
        \label{f:SaulPhi2}
        \end{center}
        \end{figure}

Finally, we define $\Phi: S \times D^2 \rightarrow M$. Suppose $x \in S$ and $p \in D^2$. If $f^i$ is non-zero at $p$ and $x \in A^i$ then we define $\Phi _p (x)=\gamma ^i (x,f^i(p))$, where $\Phi_{p}(x)=\Phi(x,p)$. Similarly, if $g^j$ is non-zero at $p$ and $x \in B^j$ then we define  $\Phi _{p}(x)=\delta ^j (x,g^j(p))$.  If $x$ is a point of $S$ not in such an $A^i$ or $B^j$ then we define $\Phi _{p}(x)=x$. 

The proof is now complete by making the following observations:
\begin{itemize}
	\item If $f^i$ and $f^k$ are non-zero at $p$ (where $k >i$) then $k=i+1$ and $g^j(p)=0$ for all $j$. Since $A^i$ and $A^{i+1}$ are disjoint the function $\Phi _{p}$ is well defined. (Ambient isotopies with disjoint supports commute.)
	\item A similar statement holds if $g^j$ and $g^l$ are non-zero at $p$. 
	\item If $f^i$ and $g^j$ are non-zero at $p$ then either $i=j=0$ or $i=n$ and $j=m$. Since $A^0 \cap B^0=A^n \cap B^m=\emptyset$ the function is again well defined. 
	\item For $p$ near the center of $D^2$ the function $\Phi _p$ is the identity on $S$. 
	\item At each $p \in \partial D^2$ at least one of the functions $\{f^i\}$ or $\{g^j\}$ is 1. Assume this is true of $f^i(p)$. If $x \in S \cap A^i$ then $\Phi _{p}(x)=\gamma ^i (x,1)$. Note that this is disjoint from the disk $V^i$. Hence, for every point $p \in \partial D^2$ the surface $S_p$ is disjoint from at least one of the disks $\{V^i\}$ or $\{W^j\}$. 
	\item At $\theta_0$ both $f^0$ and $g^0$ are 1. Hence, the surface $S _{\theta_0}$ is disjoint from both $V^0=V_0$ and $W^0=W_0$. Similarly, at $\theta _1$ both $f^n$ and $g^m$ are 1. Hence, the surface $S _{\theta_1}$ is disjoint from both $V^n=V_1$ and $W^m=W_1$.
\end{itemize}
\end{proof}

We now perturb $\Phi$ slightly so that it is smooth and in general position with respect to $H$, and again denote the new function as $\Phi$. Consider the set $\Sigma=\{x \in D^2 | S_x $ is not transverse to $H\}$. If $\Phi$ is in general position with respect to $H$, then Cerf theory (see \cite{cerf:68}) tells us that $\Sigma$ is homeomorphic to a graph, and the maximum valence of each vertex of this graph is 4. We will use these facts later.

\subsection{Labelling $D^2$.}\

A {\it region} of $D^2$ is a component of $D^2 - \Sigma$. Let $x$ be any point in some region. We will label this region from the set $\{\VV_0, \WW_0, \VV_1, \WW_1 \}$ as follows. The region containing $x$ will have the label 
\begin{itemize}
	\item $\VV_0$ if there is a loop of $H \cap S_x$  which bounds a compressing disk $V \subset \CV$ for $H$ such that $V \in C_0$. 
	\item $\WW_0$ if there is a loop of $H \cap S_x$  which bounds a compressing disk $W \subset \CW$ for $H$ such that $W \in C_0$. 
	\item $\VV_1$ if there is a loop of $H \cap S_x$  which bounds a compressing disk $V \subset \CV$ for $H$ such that $V \in C_1$. 
	\item $\WW_1$ if there is a loop of $H \cap S_x$  which bounds a compressing disk $W \subset \CW$ for $H$ such that $W \in C_1$.
\end{itemize}

\begin{clm}
\label{c:UnlabelledRegionImpliesDone}
If some region is unlabeled then the conclusion of Lemma \ref{l:EssentialIntersection} follows.
\end{clm}

\begin{proof}
If any region remains unlabeled then there is no loop of $H \cap S_x$  which bounds a compressing disk for $H$. We claim that in such a situation every loop of $H \cap S_x$ is inessential on both surfaces, and hence we can remove all intersections by a sequence of surgeries, using a standard innermost disk argument. The conclusion of Lemma \ref{l:EssentialIntersection} thus follows.

Suppose $H \cap S_x$ contains a loop $\gamma$ which is essential on $H$. As $S$ is a sphere or disk the loop $\gamma$ bounds a subdisk $D$ of $S$. Let $\gamma'$ denote a loop of $H \cap S_x$ which is innermost on $D$ among all loops that are essential on $H$ (possibly $\gamma'=\gamma$). The loop $\gamma'$ thus bounds a subdisk $D'$ of $S$ whose interior meets $H$ in a collection of loops that are inessential on both surfaces. Hence we may do a sequence of surgeries on $D'$ to obtain a compressing disk for $H$. It would follow that the region containing $x$ has a label. 
\end{proof}

\begin{clm}
\label{c:same}
No region can have both of the labels $\VV_i$ and $\WW_{1-i}$.
\end{clm}

\begin{proof}
Let $x$ be a point in a region with the labels $\VV_0$ and $\WW_1$. Let $V \subset \CV$ and $W \subset \CW$ be disks whose existence is implied by these labels. Hence, $V \in C_0$ and $W \in C_1$. But both $\partial V$ and $\partial W$ are contained in $H \cap S_x$. Thus they are either disjoint or equal (in which case they can be made disjoint). This contradicts the definition of $C_0$ and $C_1$, as then $(V,W)$ would be a reducing pair.
\end{proof}

\begin{clm}
\label{c:adjacent}
If a region has the label $\VV_i$ then no adjacent region can have the label $\WW_{1-i}$.
\end{clm}

The proof is similar to the argument of Gabai used in Claim \ref{c:Gabai}.

\begin{proof}
Suppose the region $\mathcal R_0$ has the label $\VV_0$ and $\mathcal R_1$ is an adjacent region with the label $\WW_1$. Let $x_i$ be some point in $\mathcal R_i$. Let $p:I \rightarrow D^2$ be an embedded path connecting $x_0$ to $x_1$, which does not wander into any region other than $\mathcal R_0$ or $\mathcal R_1$. The fact that $\mathcal R_0$ has the label $\VV _0$ implies that for each $t$ for which $p(t)$ is in $\mathcal R_0$ there is a compressing disk $V_t \subset \CV$ for $H$ such that $V_t \in C_0$ and $\partial V_t \subset H \cap S_{p(t)}$. Similarly, for each $t$ for which $p(t)$ is in $\mathcal R_1$ there is a compressing disk $W_t \subset \CW$ for $H$ such that $W_t \in C_1$ and $\partial W_t \subset H \cap S_{p(t)}$.

As $t$ increases from $0$ to $1$, we see a moment, $t_*$, when $S _{p(t_*)}$ does not meet $H$ transversely ({\it i.e.} $t_*$ corresponds to the place where the path $p$ crosses an edge of $\Sigma$). At $t_*$ we simultaneously see the disappearance of $\partial V_t$, and the appearance of $\partial W_t$. (Otherwise, $\mathcal R_1$ would have both the labels $\VV_0$ and $\WW_1$, which we ruled out in Claim \ref{c:same}.) We conclude that as $t$ approaches $t_*$ from below we see $\partial V_t$ become tangent to itself, or to another loop of $H \cap S_{p(t)}$. Similarly, as $t$ approaches $t_*$ from above we see $\partial W_t$ become tangent to itself, or to another loop of $H \cap S_{p(t)}$. As only one such tangency occurs for each $t$ on an edge of $\Sigma$ we see that 
\[\lim \limits _{t \rightarrow t_*^-} \partial V_t \cap \lim \limits _{t \rightarrow t_*^+} \partial W_t \ne \emptyset\]
as in Figure \ref{f:gamma}. 

        \begin{figure}[htbp]
        \psfrag{g}{$\partial V_t$}
        \psfrag{h}{$\partial W_t$}
        \psfrag{-}{$t<t_*$}
        \psfrag{+}{$t>t_*$}
        \psfrag{0}{$t=t_*$}
        \vspace{0 in}
        \begin{center}
        \includegraphics[width=3 in]{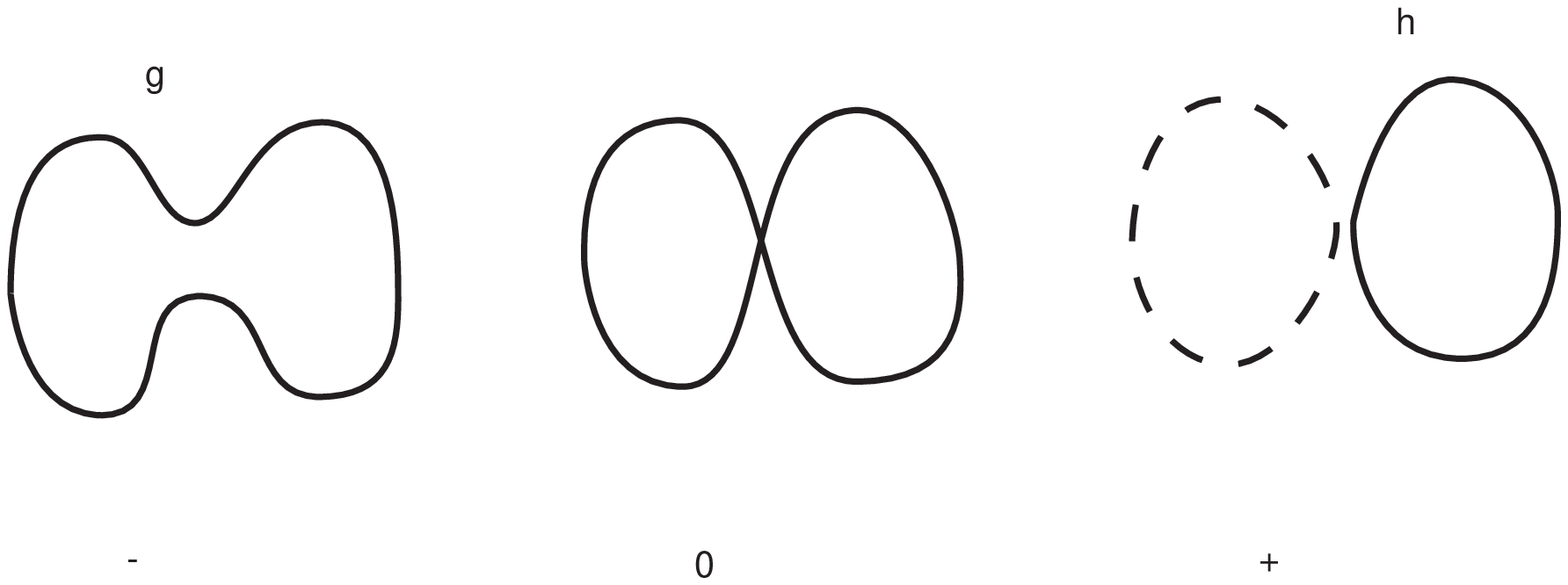}
        \caption{}
        \label{f:gamma}
        \end{center}
        \end{figure}

Since $V_t$ and $W_t$ are on opposite sides of $H$, we see from Figure \ref{f:gamma} that $\partial V_{t_*-\epsilon}$ can be made disjoint from $\partial W_{t_*+\epsilon}$ (since $H$ is orientable). This contradicts the fact that  $V_t \in C_0$ and $W_t \in C_1$.
\end{proof}

We now assume, to obtain a contradiction, that there are no unlabeled regions.

\subsection{The 2-complex $\Pi$ and a map from $D^2$ to $\Pi$.}\
                                                                        
Let $\Pi$ be the labelled 2-complex depicted in Figure \ref{f:pi}. Let $\Sigma '$ be the dual graph of $\Sigma$. Map each vertex of $\Sigma '$ to the point of $\Pi$ with the same label(s) as the region of $D^2$ in which it sits. Claim \ref{c:same} assures that this map is well defined on the vertices of $\Sigma '$. 

        \begin{figure}[htbp]
        \psfrag{1}{$\VV_0$}
        \psfrag{2}{$\VV_0$,$\VV_1$}
        \psfrag{3}{$\VV_1$}
        \psfrag{4}{$\WW_0$}
        \psfrag{5}{$\WW_0$,$\WW_1$}
        \psfrag{6}{$\WW_1$}
        \psfrag{7}{$\VV_0$,$\WW_0$}
        \psfrag{8}{$\VV_1$,$\WW_1$}
        \vspace{0 in}
        \begin{center}
        \includegraphics[width=3 in]{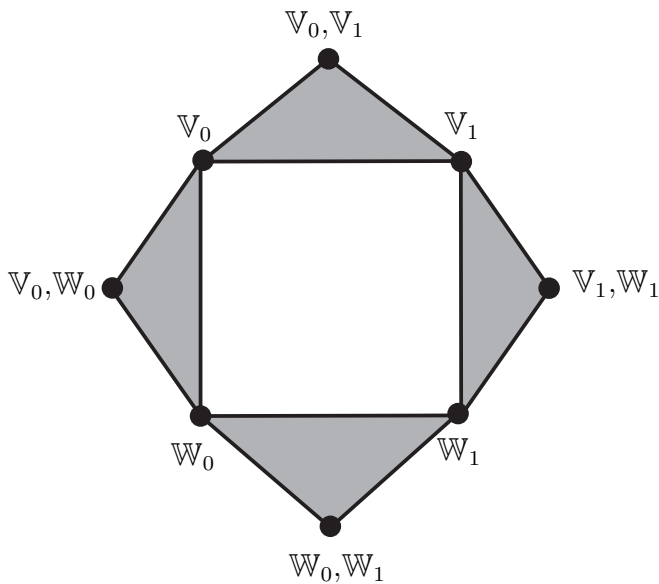}
        \caption{The 2-complex, $\Pi$.}
        \label{f:pi}
        \end{center}
        \end{figure}

Similarly, map each edge of $\Sigma '$ to the 1-simplex of $\Pi$ whose endpoints are labelled the same. Claim \ref{c:adjacent} guarantees that this, too, is well defined.

\begin{clm}
The map to $\Pi$ extends to all of $D^2$.
\end{clm}

\begin{proof}
Note that the maximum valence of a vertex of $\Sigma$ is four. Hence, the boundary of each region in the complement of $\Sigma '$ gets mapped to a 1-cycle with at most four vertices in $\Pi$. Inspection of Figure \ref{f:pi} shows that there is only one such cycle which is not null homologous. Hence, we must rule out the possibility that there are four regions around a common vertex $x_*$ of $\Sigma$, each with only one label, where all such labels are distinct. As in the proof of Claim \ref{c:adjacent} this implies that each edge of $\Sigma$ incident to $x_*$ corresponds to a saddle tangency. Hence, $S_{x_*} \cap H$ is a graph with exactly two valence four vertices and simple closed curves. 

        \begin{figure}[htbp]
        \vspace{0 in}
        \begin{center}
        \includegraphics[width=2 in]{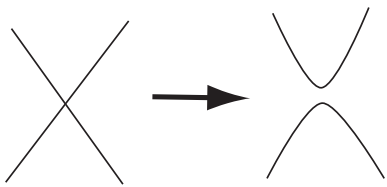}
        \caption{Resolving a vertex of $S_{x_*} \cap H$.}
        \label{f:resolve}
        \end{center}
        \end{figure}

For $x$ in the interior of a region which meets $x_*$ the set $S_x \cap H$ is obtained from $S_{x_*} \cap H$ by some resolution of its vertices (see Figure \ref{f:resolve}). There are exactly four possible ways to resolve two vertices. As there are four regions with different labels around $x_*$ we must see all four resolutions. However, the orientability of $H$ guarantees that {\it some} resolution will consist of loops  that can be made disjoint from all components of all other resolutions. There are now four symmetric cases. Suppose, for example, that such a resolution contains the boundary of a disk $V \subset \CV$ such that $V \in C_0$. We know some other resolution contains that boundary of a disk $W \subset \CW$ such that $W \in C_1$. But this contradicts the fact that $\partial V$ can be made disjoint from $\partial W$. The proof is now complete by symmetry. 
\end{proof}

\subsection{Finding an unlabelled region.}\

To obtain a contradiction to our assumption that all regions are labeled it suffices to prove that the map from $D^2$ to $\Pi$, when restricted to $\partial D^2$, induces a non-trivial map on homology. To this end we must examine the possibilities for the labels of the regions adjacent to $\partial D^2$. 

\begin{clm}
\label{c:DisjointOnBoundary}
If $S_x$ is disjoint from a compressing disk $V \subset \CV$ for $H$ such that $V \in C_i$ then the region containing $x$ does not have the label $\WW_{1-i}$. Similarly, if $S_x$ is disjoint from a compressing disk $W \subset \CW$ for $H$ such that $W \in C_i$ then the region containing $x$ does not have the label $\VV_{1-i}$.
\end{clm}

\begin{proof}
Assume $S_x$ is disjoint from a disk $V \subset \CV$, where $V \in C_i$. If the region containing $x$ has the label $\WW_{1-i}$ then there is a loop $\gamma$ of $S_x \cap H$ which bounds a compressing disk $W \subset \CW$ for $H$  such that $W \in C_{1-i}$. Then $\partial W \cap \partial V=\emptyset$. This contradicts the fact that $(V,W)$ is not a reducing pair. The proof is complete by symmetry.
\end{proof}

Recall the arcs $U$ and $L$ of $\partial D^2$ from Claim \ref{c:PhiExists}. For any point $x$ near $U$ the surface $S_x$ is disjoint from some compressing disk for $H$ in $\CV$. Similarly, for any $x$ near $L$ the surface $S_x$ is disjoint from a compressing disk for $H$ in $\CW$.

\begin{clm}
\label{c:AdjacentOnBoundary}
Suppose $\mathcal R_0$ and $\mathcal R_1$ are regions such that $\mathcal R_0 \cap U$ is adjacent to $\mathcal R_1 \cap U$. If $\mathcal R_0$ has the label $\WW_0$ then $\mathcal R_1$ cannot have the label $\mathcal \WW_1$. Similarly, if $\mathcal R_0 \cap L$ is adjacent to $\mathcal R_1 \cap L$ and  $\mathcal R_0$ has the label $\VV_0$ then $\mathcal R_1$ cannot have the label $\VV_1$.
\end{clm}

\begin{proof}
Suppose $\mathcal R_0 \cap U$ is adjacent to $\mathcal R_1 \cap U$, the label of $\mathcal R_0$ is $\WW_0$, and the label of $\mathcal R_1$ is $\WW_1$. Let $p$ be the point $\mathcal R_0 \cap \mathcal R_1 \cap U$. By Claim \ref{c:PhiExists} the surface $S_p$ is disjoint from a compressing disk $V \subset \CV$ for $H$. Hence, for all points $x$ near $p$ the surface $S_x$ will be disjoint from $V$. But every neighborhood of $p$ contains points in both $\mathcal R_0$ and $\mathcal R_1$. For $x$ near $p$ and in $\mathcal R_0$ Claim \ref{c:DisjointOnBoundary} implies $V \in C_0$. But if $x \in \mathcal R_1$ then Claim \ref{c:DisjointOnBoundary} implies $V \in C_1$. As $C_0$ and $C_1$ partition the compressing disks for $H$ the disk $V$ cannot be in both.
\end{proof}

\begin{figure}[htbp]
        \psfrag{1}{$\VV_0$}
        \psfrag{2}{$\VV_0$,$\VV_1$}
        \psfrag{3}{$\VV_1$}
        \psfrag{4}{$\WW_0$}
        \psfrag{5}{$\WW_0$,$\WW_1$}
        \psfrag{6}{$\WW_1$}
        \psfrag{7}{$\VV_0$,$\WW_0$}
        \psfrag{8}{$\VV_1$,$\WW_1$}
        \psfrag{a}{$\theta_0$}
        \psfrag{b}{$\theta_1$}
        \psfrag{U}{$U$}
        \psfrag{L}{$L$}
        \vspace{0 in}
        \begin{center}
        \includegraphics[width=4 in]{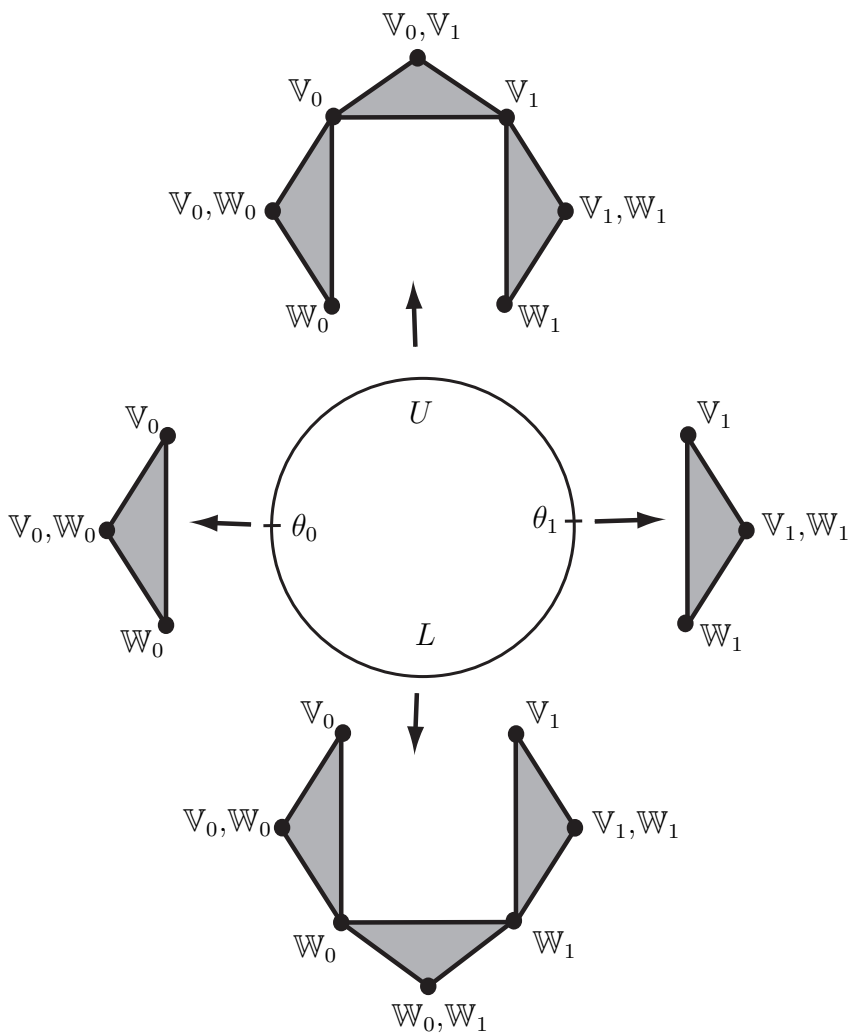}
        \caption{The map from $\partial D^2$ to $\Pi$ is non-trivial on homology.}
        \label{f:rpi}
        \end{center}
        \end{figure}

We now examine the properties of $\Phi$ listed in Claim \ref{c:PhiExists} and the complex $\Pi$. If $x$ is in a region containing the point $\theta_0$ then $S_x$ is disjoint from $V_0$ and $W_0$. It follows from Claim \ref{c:DisjointOnBoundary} that this region can not have either of the labels $\VV_1$ or $\WW_1$, and hence must get mapped to the left triangle of $\Pi$. Similarly, a region containing the point $\theta_1$ gets mapped to the right triangle of $\Pi$. If $x$ is in a region containing a point of $U$ then $S_x$ is disjoint from some compressing disk for $H$ in $\CV$. It follows from Claim \ref{c:DisjointOnBoundary} that such a region cannot get both of the labels $\WW_0$ and $\WW_1$, and hence does not get mapped to the bottom point of $\Pi$. It follows from Claim \ref{c:AdjacentOnBoundary} that no adjacent pair of regions next to $U$  gets mapped to the horizontal edge of the bottom triangle of $\Pi$. We conclude that the arc $U$ gets mapped to the union of the left, right, and top triangles of $\Pi$. A symmetric argument shows that $L$ gets mapped to the union of the left, right, and bottom triangles of $\Pi$. We summarize these observations in Figure \ref{f:rpi}. It follows that the map from $D^2$ to $\Pi$, when restricted to $\partial D^2$, is non-trivial on homology, a contradiction. We conclude that there must have been an unlabeled region, as this was the only assumption we made in constructing the map. Lemma \ref{l:EssentialIntersection} now follows from Claim \ref{c:UnlabelledRegionImpliesDone}.
\end{proof}

\section{Generalized Heegaard Splittings}
\label{s:GHSdefinition}

Theorem \ref{t:main} is essentially proved by a complex sequence of handle slides. However, we find that a handle structure is more cumbersome to deal with than a {\it generalized Heegaard splitting}, defined presently. The relationship between handle structures and generalized Heegaard splittings is made explicit in Example \ref{ex:GHSvsHandles}. 

\begin{dfn}
\label{d:GHS}
A {\it generalized Heegaard splitting (GHS)} $H$ of a 3-manifold $M$ is a pair of sets of pairwise disjoint, connected surfaces,  $\thick{H}$ and $\thin{H}$ (called the {\it thick levels} and {\it thin levels}, respectively), which satisfy the following conditions. 
	\begin{enumerate}
		\item Each component $M'$ of $M-\thin{H}$ meets a unique element $H_+$ of $\thick{H}$, and $H_+$ is a Heegaard surface in $M'$. Henceforth we will denote the closure of the component of $M-\thin{H}$ that contains an element $H_+ \in \thick{H}$ as $M(H_+)$. 
		\item As each Heegaard surface $H_+ \subset M(H_+)$ is oriented, say, by an outward pointing normal vector, we can consistently talk about the points of $M(H_+)$ that are ``above" $H_+$ or ``below" $H_+$. Suppose $H_-\in \thin{H}$. Let $M(H_+)$ and $M(H_+')$ be the submanifolds on each side of $H_-$ (a priori it is possible that $M(H_+)=M(H_+')$). Then $H_-$ is below $H_+$ if and only if it is above $H_+'$.
		\item There is a partial ordering on the elements of $\thin{H}$ which satisfies the following: Suppose $H_+$ is an element of $\thick{H}$, $H_-$ is a component of $\partial M(H_+)$ above $H_+$, and $H'_-$ is a component of $\partial M(H_+)$ below $H_+$. Then $H_- > H'_-$.
	\end{enumerate}
\end{dfn}

%What about parallel thick and thin levels?

\begin{ex}
Suppose $H_+$ is a Heegaard surface in a closed 3-manifold $M$. Then a GHS $H$ of $M$ is given by $\thick{H}=\{H_+\}$ and $\thin{H}=\emptyset$.
\end{ex}

\begin{figure}[htbp]
        \psfrag{F}{$F_{1}$}
        \psfrag{G}{$F_{2}$}
        \psfrag{A}{$H_{+}^{1}$}
        \psfrag{B}{$H_{+}^{2}$}
        \vspace{0 in}
        \begin{center}
        \includegraphics[width=2 in]{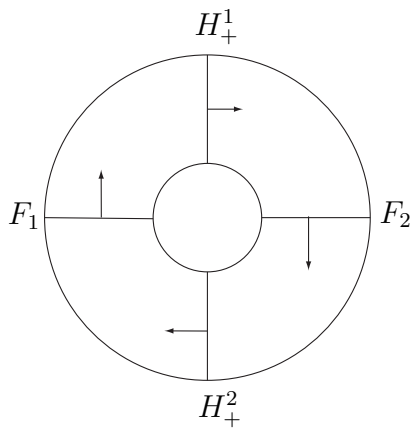}
        \caption{The orientations on $H_{+}^{1}$ and $H_{+}^{2}$ are given by normal vectors. If the orientations on $F_{1}$ and $F_{2}$ are chosen to be consistent with condition (2) of Definition \ref{d:GHS} (as pictured) then they will be inconsistent with condition (3).}
        \label{f:OrientedGHS}
        \end{center}
        \end{figure}

\begin{ex}
Suppose $M$ is a 3-manifold, $H_+$ is a Heegaard surface in $M$, and $\partial M$ has two homeomorphic components separated by $H_+$. Let $M^1$ and $M^2$ denote two homeomorphic copies of $M$, and $H^i_+$ the image of $H_+$ in $M^i$. Let $N$ denote the connected 3-manifold obtained from $M^1$ and $M^2$ by identifying the boundary components in pairs by some homeomorphisms. The image of the boundary components in $N$ are $F_1$ and $F_2$. We may attempt to define a GHS $H$ by setting $\thick{H}=\{H^1_+, H^2_+\}$ and $\thin{H}=\{F_1,F_2\}$. But this $H$ is NOT a GHS of $N$, since there is no way to consistently satisfy conditions (2) and (3) of the definition. If orientations are chosen on $H_+^1$ and $H_+^2$ to satisfy condition (2), then the relation specified by condition (3) will imply both $F_1>F_2$ and $F_2>F_1$. See Figure \ref{f:OrientedGHS}. 
\end{ex}

\begin{figure}[htbp]
        \psfrag{0}{0- and 1-handles}
        \psfrag{2}{2- and 3-handles}
        \psfrag{T}{Thick levels}
        \psfrag{t}{Thin levels}
        \vspace{0 in}
        \begin{center}
       \includegraphics[width=3.75 in]{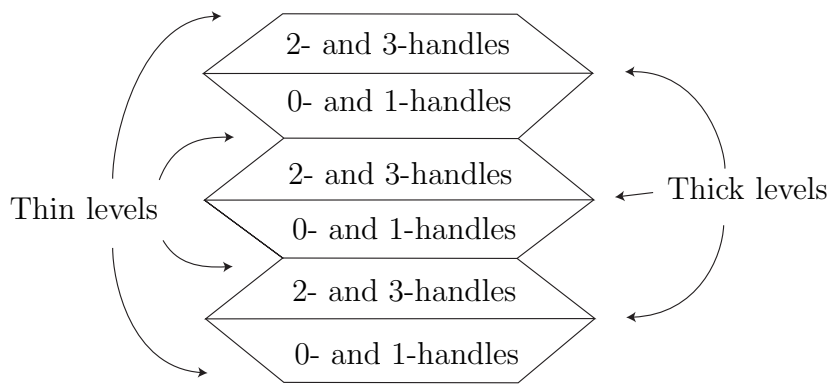}
       \caption{The GHS associated with a handle structure.}
        \label{f:handles}
        \end{center}
        \end{figure}

\begin{ex}
\label{ex:GHSvsHandles}
The original formulation of a GHS was given by Scharlemann and Thompson in \cite{st:94}. They defined these structures as being ``dual" to handle structures. To be explicit, assume $M$ is built by starting with  0-handles and attaching 1-handles, 2-handles, and 3-handles in any order. The set of thick and thin levels of the GHS associated with this handle structure appears at the various interfaces of the 1- and 2-handles, as indicated in Figure \ref{f:handles}. Note that some of the horizontal lines in the figure may represent a disconnected surface. Each component of such a surface will be either a thin or thick level.  
\end{ex}

\begin{dfn}
Suppose $H$ is a GHS of a 3-manifold $M$ with no 3-sphere components. Then $H$ is {\it strongly irreducible} if each element $H_+ \in \thick{H}$ is strongly irreducible 
%or incompressible???
in $M(H_+)$.  We say $H$ is {\it critical} if there is a unique element $H_*\in \thick{H}$ which is critical in $M(H_*)$, and every other element $H_+ \in \thick{H}$ is strongly irreducible in $M(H_+)$.
\end{dfn}

%What about products?? What if $M\cong F \times S^1$, and $\thick{H}=F \times\{p\}$, $\thin{H}=F \times \{q\}$?

\begin{lem}
\label{l:IncompressibleThinLevels}
Suppose $H$ is a strongly irreducible or critical GHS of $M$. Then each element of $\thin{H}$ is  incompressible in $M$.
\end{lem}

The strongly irreducible case is motivated by \cite{st:94}.

\begin{proof}
Let $H$ be a strongly irreducible or critical GHS. Choose a compressing disk $D$ for some thin level whose interior meets the union of all thin levels a minimal number of times. Let $\alpha$ denote an intersection loop of $D$ with the set of thin levels which is innermost on $D$ (possibly $\alpha=\partial D$). Let $H_-$ denote the thin level that contains $\alpha$. Suppose first that  $\alpha$ is inessential on $H_-$, bounding a subdisk $A$ of this surface. Let  $B$ be a subdisk of $A$ bounded by a loop of $D \cap A$ which is innermost on $A$. Then we may use $A$ to surger $D$, producing a new compressing disk for some thin level which meets the union of all thin levels fewer times. 

Now suppose $\alpha$ is essential on $H_-$. The loop $\alpha$ bounds a subdisk $A'$ of $D$ which lies in $M(H_+)$, for some $H_+ \in \thick{H}$. By assumption $H_+$ is either a strongly irreducible or a critical Heegaard surface in $M(H_+)$. 

By Lemma \ref{l:EssentialIntersection} there is a disk in $M(H_+)$, with the same boundary as $A'$, which misses $H_+$. This disk is thus a compressing disk for the negative boundary of a compression body, a contradiction. 
\end{proof}

\begin{lem}
\label{l:SIorCritImpliesIrreducible}
If $M$ admits a strongly irreducible or critical GHS then $M$ is irreducible. 
\end{lem}

\begin{proof}
Let $H$ denote a strongly irreducible or critical GHS. By Lemma \ref{l:IncompressibleThinLevels} the union of the set of thin levels of $H$ forms a (disconnected) incompressible surface in $M$. By Lemma \ref{l:EssentialIntersection} there is thus an essential sphere $S'$ which is disjoint from $\thin{H}$. Hence $S' \subset M(H_+)$ for some $H_+ \in \thick{H}$. Again by Lemma \ref{l:EssentialIntersection} there must be an essential sphere $S'' \subset M(H_+)$ that misses $H_+$. But then $S''$ is an essential sphere in a compression body, which can not exist. 
\end{proof}

\section{Reducing GHSs}
\label{s:WRdefinition}

\begin{figure}[htbp]
        \psfrag{1}{$G_+/D$}
        \psfrag{2}{$G_+/E$}
        \psfrag{3}{$G_+/DE$}
        \psfrag{G}{$G_+$}
        \psfrag{E}{$E$}
        \psfrag{D}{$D$}
        \vspace{0 in}
        \begin{center}
       \includegraphics[width=3.5 in]{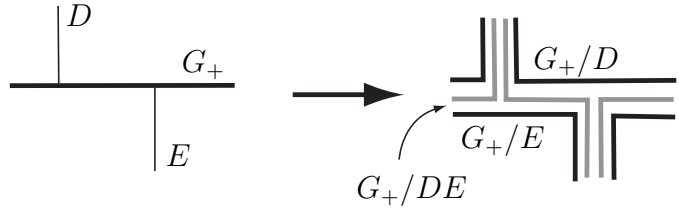}
       \caption{The first step in defining a weak reduction.}
        \label{f:WeakReduction}
        \end{center}
        \end{figure}

We now define a way to take a GHS $G$ of a 3-manifold $M$ which is not strongly irreducible and obtain a ``simpler" GHS $H$. The new GHS will be of a manifold $M'$ which is obtained from $M$ by cutting along a sphere and capping off the resulting boundary components with a pair of 3-balls. This procedure will be called {\it weak reduction}, and can be performed whenever there is a thick level $G_+ \in \thick{G}$ and a reducing pair $(D,E)$ for $G_+$ in $M(G_+)$. We assume $M$ has no 3-sphere components. We define the new GHS by the following algorithm:

\begin{enumerate}
	\item Initially set  
	\[\thick{H}=\thick{G} -\{G_+\} \cup \{G_+ / D, G_+ / E\},\]
	\[\thin{H} =\thin{G} \cup \{G_+/DE\},\ \mbox{and}\ M'=M.\]
	See Figure \ref{f:WeakReduction}.
	\item If there is now a sphere $S \in \thin{H}$ then cut $M'$ along $S$, cap off the resulting sphere boundary components with 3-balls, and remove $S$ from $\thin{H}$. 
	\item If $M'$ now contains a 3-sphere component then delete it, along with all elements of $\thick{H}$ and $\thin{H}$ that lie in this component. 
	\item If there are elements $H_+ \in \thick{H}$ and $H_- \in \thin{H}$ that cobound a product region $P$ of $M'$ such that $P \cap \thick{H}=H_+$ and $P \cap \thin{H}=H_-$ then remove $H_+$ from $\thick{H}$ and $H_-$ from $\thin{H}$. 
\end{enumerate}

The first step of the above algorithm is illustrated in Figure \ref{f:WeakReduction}.	

%\begin{lem}
%Weak reduction produces a GHS.
%\end{lem}

\begin{ex}
\label{e:destab1}
Suppose $\thick{H}$ contains a single element $H_+$ and $\thin{H}=\emptyset$, so that $H_+$ is a Heegaard surface for $M$. Suppose $(D,E)$ is a weak reduction for $H$ and $\partial D$ and $\partial E$ cobound an annulus $A$ of $H_+$. If the sphere $D \cup A \cup E$ is separating in $M$ then the inverse of the weak reduction is a connected sum, in the sense of Definition \ref{d:ConnectedSum}. 

Suppose $H_+$ is obtained from some lower genus Heegaard surface by a stabilization. Recall that a stabilization is a connected sum with the standard genus one Heegaard surface in $S^3$. Hence there will be a weak reduction for $H$ which undoes this.
\end{ex}

In the previous example we saw that a sphere may appear at a thin level in the process of weak reduction when $\partial D$ and $\partial E$ were parallel. The next example shows that there is another way for spheres to crop up at thin levels during weak reduction.

\begin{ex}
\label{e:destab2}
Let $M_1$ and $M_2$ be 3-manifolds with genus one Heegaard splittings (possibly one or both are $S^3$). Let $M=M_1\#M_2$. Then, as in Example \ref{e:LensSpaces}, there is a genus two Heegaard splitting $H_+$ of $M$ which is the connected sum of genus one splittings of $M_1$ and $M_2$. There are compressing disks $D \subset M_1$ and $E \subset M_2$ on opposite sides of $H_+$ in $M$. Hence the pair $(D,E)$ is a reducing pair. In forming the weak reduction one of the first steps is to initially add $H_+/DE$ to $\thin{H}$. But this surface is a sphere, which will get removed in subsequent steps. 

If $M_2 \cong S^3$ then $H_+$ was a stabilization of the splitting of $M_1$ that we started with. The weak reduction given by $(D,E)$  has undone this stabilization. 
\end{ex}

The preceding pair of examples motivates us to make the following definition.

\begin{dfn}
The weak reduction of a GHS given by the reducing pair $(D,E)$ for the thick level $H_+$ is called a {\it destabilization} if $H_+/DE$ contains a sphere which bounds a ball. 
\end{dfn} 

There are two ways that the weak reduction given by $(D,E)$ can be a destabilization. The first is when $\partial D$ and $\partial E$ are parallel on $H_+$, as in Example \ref{e:destab1}. The second is when $\partial D$ and $\partial E$ are non-parallel but $H_+$ has genus two, as in Example \ref{e:destab2}.

\begin{ex}
\label{e:S^2xS^1}
Consider $S^2 \times S^1$. If $l$ is a loop in $S^2$ then $l \times S^1$ is a Heegaard torus, $H$. The loop $l$ bounds disks $D$ and $E$ on opposite sides of $H$, which can be isotoped to be disjoint. In the process of weak reduction the manifold will get cut along an essential 2-sphere, and the resulting boundary components will get capped off by 3-balls. The result is $S^3$, which will then get deleted. Hence, weak reducing along the pair $(D,E)$ produces the empty GHS of the empty set.
\end{ex}

We now define a partial ordering of GHSs. Any partial ordering will suffice, as long as a weak reduction produces a smaller GHS, and any monotonically decreasing sequence of GHSs must terminate. This is motivated by \cite{st:94}.

\begin{dfn}
\label{d:GHScomplexity}
If $H$ is a GHS then let $c(H)$ denote the set of genera of the elements of $\thick{H}$, where repeated integers are included, and the set is put in non-increasing order. We compare two such sets lexicographically. If $H^1=\{H^1_i\}$ and $H^2=\{H^2_j\}$ are two GHSs then we say $H^1<H^2$ if $c(H^1)<c(H^2)$. 
\end{dfn}

\begin{lem}
\label{l:WRimpliesDecrease}
Suppose $H^1$ and $H^2$ are GHSs and $H^1$ is obtained from $H^2$ by a weak reduction. Then $H^1 < H^2$. 
\end{lem}

\begin{proof}
In the definition of weak reduction we remove $H_+$ from $\thick{H}$ and replace it with (possibly multiple) surfaces of smaller genus. Some of these surfaces may then appear at the boundary of a product submanifold or in an $S^3$ component, and so will be removed. In any case, $c(H)$ has gone down under the lexicographical ordering.
\end{proof}

\section{Swapping Weak Reductions}

There will be many times when we will want to alter a sequence of weak reductions. In this section we list a few basic ways to do this. 

\begin{lem}
\label{l:WRswapSeparating}
Suppose $(D,E)$ and $(D',E)$ are different weak reductions for a GHS $G$. Let $G_+$ denote the thick level of $G$ which contains $\partial D$, $\partial D'$, and $\partial E$, and suppose $\partial D$ separates $\partial E$ from $\partial D'$ on $G_+$. Then performing the weak reduction $(D,E)$ on $G$ yields the same GHS as performing $(D',E)$ followed by $(D,E)$.
\end{lem}

\begin{proof}
The proof is indicated in Figure \ref{f:Separating}. After performing the weak reduction $(D,E)$ the components of the surface $G_+/D$ are thick levels. (This is illustrated in the left side of the figure.) But $\partial D'$ and $\partial E$ now lie on different components, so we cannot follow with the weak reduction $(D',E)$. 

Since $\partial D$ separates $\partial D'$ from $\partial E$ on $G_+$ it follows that $\partial D$ and $\partial E$ lie on the same component of $G_+/D'$. Hence we may follow the weak reduction $(D',E)$ with the weak reduction $(D,E)$. (This is indicated in the right side of the figure.) After removing parallel thick and thin levels  (as indicated by the shaded regions in the figure) the two GHSs that we obtain are the same. 

        \begin{figure}[htbp]
        \psfrag{D}{$D$}
        \psfrag{d}{$D'$}
        \psfrag{E}{$E$}
        \vspace{0 in}
        \begin{center}
       \includegraphics[width=4 in]{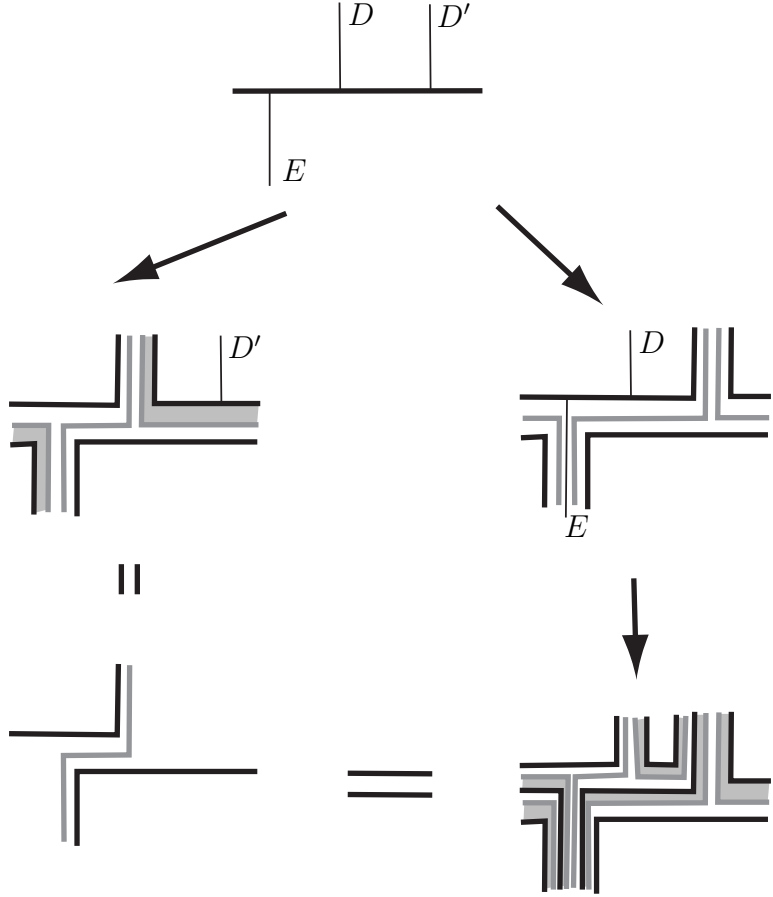}
       \caption{If $\partial D$ separates $\partial D'$ from $\partial E$ then performing the weak reduction $(D,E)$ yields the same GHS as performing $(D',E)$ followed by $(D,E)$.}
        \label{f:Separating}
        \end{center}
        \end{figure}
\end{proof}

\begin{lem}
\label{l:WRswapNonSeparating}
Suppose $(D,E)$ and $(D',E)$ are different weak reductions for a GHS $G$. Let $G_+$ denote the thick level of $G$ which contains $\partial D$, $\partial D'$, and $\partial E$, and suppose neither $\partial D$ nor $\partial D'$ separates the other from $\partial E$ on $G_+$. If, furthermore, neither $\partial D$ nor $\partial D'$ are parallel to $\partial E$ then performing the weak reduction $(D,E)$ followed by $(D',E)$ yields the same GHS as performing $(D',E)$ followed by $(D,E)$.
\end{lem}

\begin{proof}
The proof is indicated in Figure \ref{f:NonSepNonParallel}. After performing the weak reduction $(D,E)$ the curves $\partial D'$ and $\partial E$ lie on same component of $G_+/D$ so we can follow with the weak reduction $(D',E)$ (as indicated in the left side of the figure). Similarly, we may follow the weak reduction $(D'E)$ the weak reduction with $(D,E)$ (as indicated in the right side of the figure). After removing parallel thick and thin levels the two GHSs are the same. 
        \begin{figure}[htbp]
        \psfrag{D}{$D$}
        \psfrag{d}{$D'$}
        \psfrag{E}{$E$}
        \vspace{0 in}
        \begin{center}
       \includegraphics[width=4 in]{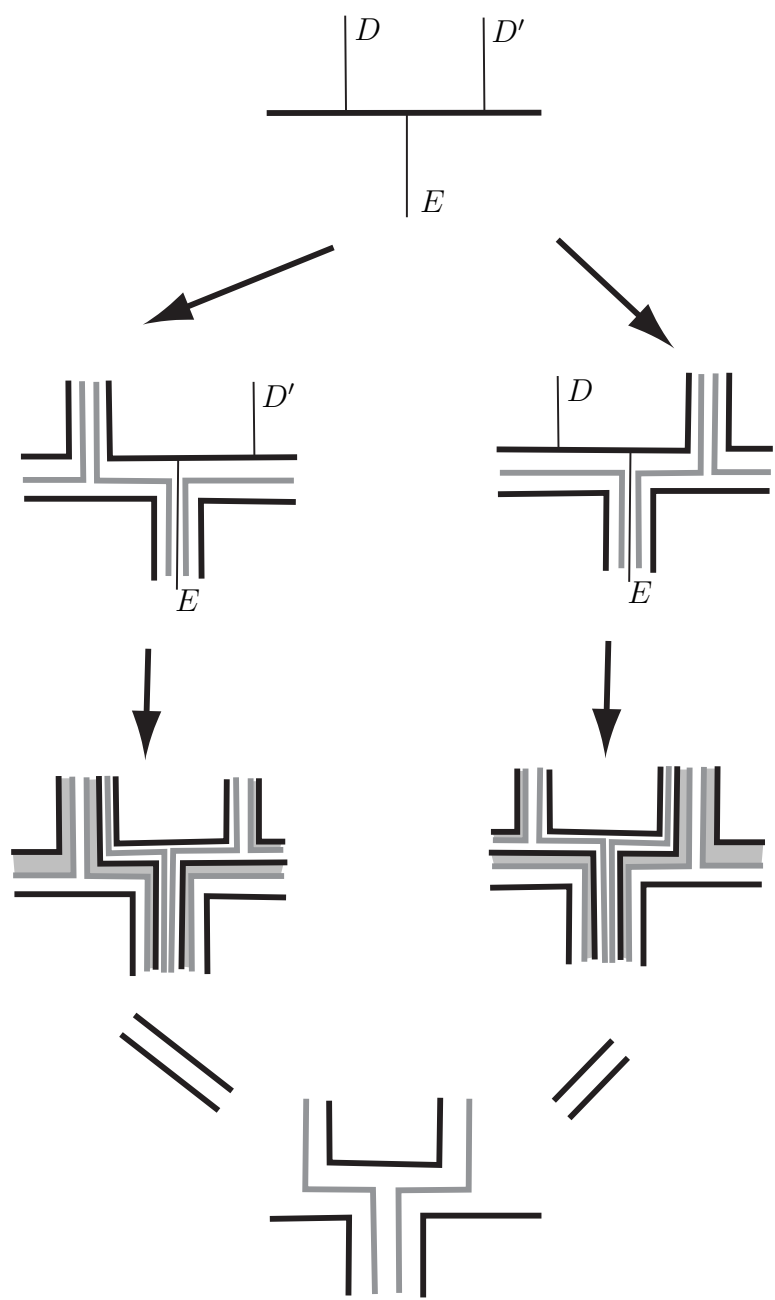}
       \caption{Performing $(D,E)$ followed by $(D',E)$ yields the same GHS as $(D',E)$ followed by $(D,E)$.}
        \label{f:NonSepNonParallel}
        \end{center}
        \end{figure}
\end{proof}

\begin{lem}
\label{l:WRswapParallel}
Suppose $(D,E)$ and $(D',E)$ are different weak reductions for a GHS $G$. Let $G_+$ denote the thick level of $G$ which contains $\partial D$, $\partial D'$, and $\partial E$. If $\partial D'$ is parallel to $\partial E$ then performing the weak reduction $(D,E)$ followed by $(D',E)$ yields the same GHS as performing $(D',E)$.
\end{lem}

\begin{proof}
The proof is indicated in Figure \ref{f:NonSepParallel}. On the left side of this figure we illustrate the following steps: (1) Performing the weak reduction $(D,E)$. After this the curves $\partial D'$ and $\partial E$ lie on same component of $G_+/D$ so we can follow with (2) the weak reduction $(D',E)$. There are now parallel thick and thin levels that can be removed (3). There is also a sphere thin level, so to complete the formation of the GHS we must (4) cut along the sphere and cap off with 3-balls. When we do this some thin level becomes parallel to a thick level (5 and 6), and so they both get removed (e). 

On the right side we (a) do $(D',E)$ first. We can not follow with $(D,E)$ because $E$ is not a compressing disk for $G_+/D'$. However, there is a sphere thin level, so we must again (b) cut and cap off with 3-balls. Once again there is now a thin level which becomes parallel to a thick level (c and d), and so they both get removed (e). 

        \begin{figure}[htbp]
        \psfrag{D}{$D$}
        \psfrag{d}{$D'$}
        \psfrag{E}{$E$}
        \psfrag{1}{$1$}
        \psfrag{2}{$2$}
        \psfrag{3}{$3$}
        \psfrag{4}{$4$}
        \psfrag{5}{$5$}
        \psfrag{6}{$6$}
        \psfrag{a}{a}
        \psfrag{b}{b}
        \psfrag{c}{c}
        \psfrag{e}{d}
        \psfrag{f}{e}
        \vspace{0 in}
        \begin{center}
       \includegraphics[width=5 in]{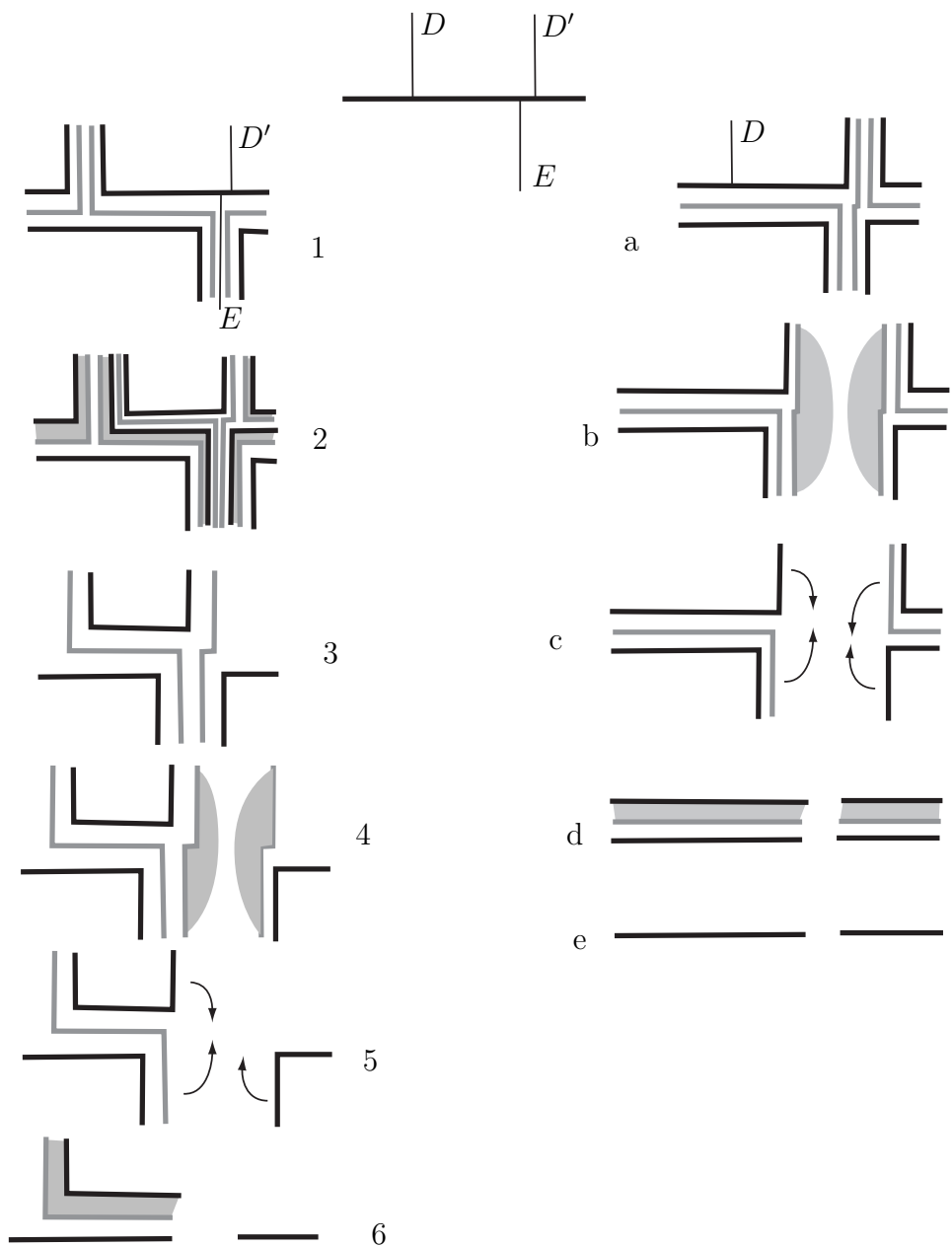}
       \caption{Performing $(D,E)$ followed by $(D',E)$ (left side) yields the same GHS as performing just $(D',E)$ (right side).}
        \label{f:NonSepParallel}
        \end{center}
        \end{figure}
\end{proof}

\section{Amalgamations}

Let $H$ be a GHS of a 3-manifold $M$. In this section we use $H$ to produce a graph $\Sigma$ in $M$. Each component of $M$ will contain exactly one component of $\Sigma$. Each component of $\Sigma$ will be the spine of a Heegaard splitting of the component of $M$ that contains it. We call the disjoint union of these Heegaard splittings the {\it amalgamation} of $H$. First, we must introduce some new notation. 

\begin{dfn}
Suppose $H$ is a GHS of $M$ and $H_{+}\in \thick{H}$. Recall that $H_{+}$ is oriented, so that we may consistently talk about those points of $M(H_+)$ that are ``above" $H_+$ and those points that are ``below." The surface $H_+$ divides $M(H_+)$ into two compression bodies. Henceforth we will denote these compression bodies as $\Wup(H_+)$ and $\Wdown(H_+)$, where $\Wup(H_+)$ is above $H_+$ and $\Wdown(H_+)$ is below. When we wish to make reference to an arbitrary compression body which lies above or below some thick level we will use the notation $\Wup$ and $\Wdown$. Define $\bdyup M(H_+)$ to be $\partial _- \Wup(H_+)$ and $\bdydown M(H_+)$ to be $\partial _- \Wdown (H_+)$. That is, $\bdyup M(H_+)$ and $\bdydown M(H_+)$ are the boundary components of $M(H_+)$ that are above and below $H_+$, respectively. 
\end{dfn} 

\begin{figure}[htbp]
        \psfrag{h}{$H_+$}
        \psfrag{W}{$\Wdown(H_+)$}
        \psfrag{w}{$\Wup(H_+)$}
        \psfrag{S}{$\Sigma$}
        \vspace{0 in}
        \begin{center}
       \includegraphics[width=3 in]{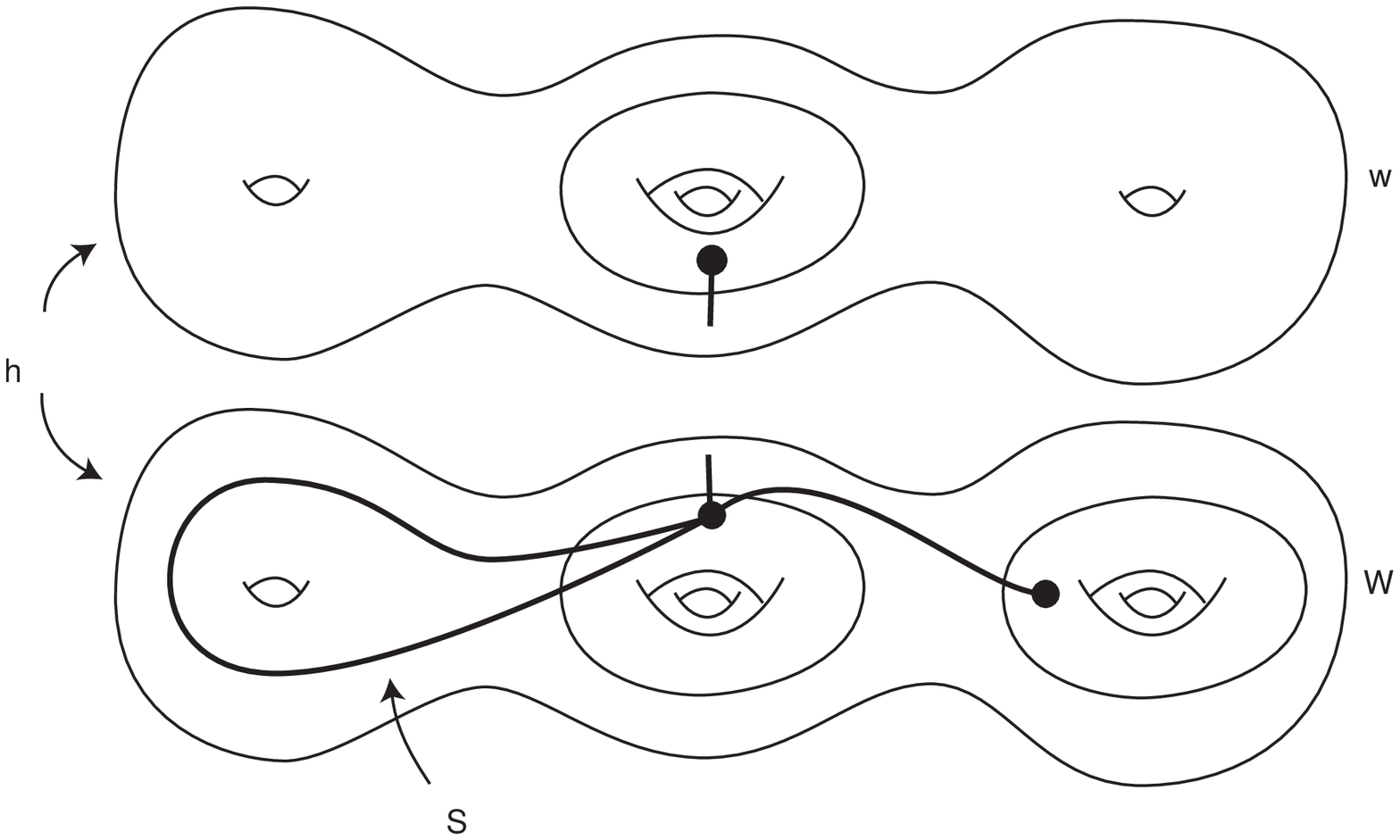}
       \caption{The intersection of $\Sigma$ with $\Wup(H_+)$ and $\Wdown(H_+)$.}
        \label{f:Amalgam}
        \end{center}
        \end{figure}

We now inductively build $\Sigma$. The intersection of $\Sigma$ with some $M(H_+)$ is depicted in Figure \ref{f:Amalgam}. First, we define a sequence of manifolds $\{M_i\}$ where
\[M_0 \subset M_1 \subset ... \subset M_n=M.\]
The submanifold $M_0$ is defined to be the disjoint union of all manifolds of the form $M(H_+)$, such that $\bdydown M(H_+)=\emptyset$. The fact that $M$ is closed and the thin levels of $H$ are partially ordered guarantees $M_0 \ne \emptyset$. Now, for each $i$ we define $M_i$ to be the union of $M_{i-1}$ and all manifolds $M(H_+)$ such that $\bdydown M(H_+) \subset \partial M_{i-1}$. Again, it follows from the partial ordering of thin levels that for some $i$ the manifold $M_i=M$. 

We now define a sequence of graphs $\Sigma_i$ in $M$. The final element of this sequence will be the desired graph $\Sigma$. 

Each $\Wdown \subset M_0$ is a handle-body. Choose a spine of each, and let $\Sigma'_0$ denote the union of these spines. The complement of $\Sigma'_0$ in $M_0$ is a (disconnected) compression body, homeomorphic to the union of the compression bodies $\Wup \subset M_0$. Now let $\Sigma_0$ be the union of $\Sigma'_0$ and one vertical arc for each component $H_-$ of $\partial M_0$,  connecting $H_-$ to $\Sigma' _0$. 

We now assume $\Sigma _{i-1}$ has been constructed and we construct $\Sigma_i$. Let $M_i'=\overline{M_i - M_{i-1}}$. For each compression body $\Wdown \subset M_i'$ choose a set of arcs $\Gamma \subset\Wdown$ such that $\partial \Gamma \subset \Sigma _{i-1} \cap \partial M_{i-1}$, and such that the complement of $\Gamma$ in $\Wdown$ is a product. Let $\Sigma' _i$ be the union of $\Sigma _{i-1}$ with all such arcs $\Gamma$. Now let $\Sigma_i$ be the union of $\Sigma'_i$ and one vertical arc for each component $H_-$ of $\partial M_i$, connecting $H_-$ to $\Sigma' _i$.

%The graph $\Sigma$ will presently be defined by its intersection with $M(H_+)$, for each $H_+ \in \thick{H}$. Each such $H_+$ divides $M(H_+)$ into compression bodies $\Wup(H_+)$ and $\Wdown(H_+)$, where $\Wup(H_+)$ is above $H_+$ and $\Wdown(H_+)$ is below. For each $H_-$ choose a point $p(H_-) \in H_-$. For each $H_+$ define $\Sigma(H_+)=\Sigma \cap M(H_+)$ as follows (see Figure \ref{f:Amalgam}):
%\begin{enumerate}
%	\item If $\partial _- \Wup(H_+) =\emptyset$ then initially let $\Sigma(H_+)$ be a spine  of $\Wup(H_+)$. If $\partial _- \Wup(H_+) \ne \emptyset$ then initially let $\Sigma(H_+)$ be a set of arcs  in $\Wup(H_+)$ that meet $\partial _-\Wup(H_+)$ in $\{p(H_-)|H_- \subset \partial _- \Wup(H_+)\}$ such that the complement of $\partial _-\Wup(H_+) \cup \Sigma(H_+)$ in $\Wup(H_+)$ is a product. 
%	\item The complement of $\partial _-\Wup(H_+) \cup \Sigma(H_+)$ in $M(H_+)$ is a compression body. Hence, for each point $p(H_-)$, where $H_- \subset \partial _-\Wdown(H_+)$ we may choose a vertical arc in $M(H_+)$ which connects $p(H_-)$ to $\partial _-\Wup(H_+) \cup \Sigma(H_+)$. If one endpoint of this arc ends up on $H_- \subset \partial _-\Wup(H_+)$ then slide it until it coincides with $p(H_-)$. 
%\end{enumerate}

\begin{lem}
If $H$ is a GHS of $M$ then each component of the graph $\Sigma$ defined above is the spine of a Heegaard splitting of the component of $M$ that contains it.
\end{lem}

\begin{proof}
Recall the sequence of manifolds $\{M_i\}$ above. We prove the lemma by showing that if the complement of $N(\Sigma)$ in $M_{i-1}$ is a union of compression bodies then the complement of $N(\Sigma)$ in $M_{i}$ is a union of compression bodies. For this it is convenient to set $M_{-1}=\emptyset$. 

%Let $N(\Sigma)$ denote a regular neighborhood of $\Sigma$ in $M$. We first show that the complement of $N(\Sigma)$ in $M_0$ is a union of handle-bodies. Suppose $\Wdown \subset M_0$ then the complement of $N(\Sigma)$ in $\Wdown$ is homeomorphic to $\{\mbox{punctured surface}\} \times I$, which is a handle-body. Now suppose $\Wup \subset M_0$. The complement of $N(\Sigma)$ in $\Wup(H_+)$ is of the form $\{\mbox{punctured surface}\} \times I$, together with a collection of 1-handles. Gluing this to the $\{\mbox{punctured surface}\} \times I$ that is the complement of $N(\Sigma)$ in $\Wdown$ thus produces a larger handle-body. 

%We now show that if the complement of $N(\Sigma)$ in $M_{i-1}$ is a union of compression bodies then the complement of $N(\Sigma)$ in $M_{i}$ is a union of compression bodies, thereby completing the proof. 

Suppose $M(H_+) \subset M_{i}'=\overline{M_{i} - M_{i-1}}$. The complement of $N(\Sigma)$ in $\Wdown(H_+)$  is a product of the form $\{\mbox{punctured surface}\} \times I$, so when we glue this on to the complement of $N(\Sigma)$ in $M_{i-1}$ we do not change the homeomorphism type of the manifold. The complement of $N(\Sigma)$ in $\Wup(H_+)$ is of the form $\{\mbox{punctured surface}\} \times I$, together with a collection of 1-handles. Gluing this to the complement of $N(\Sigma)$ in $M_{i-1} \cup \Wdown(H_+)$ thus produces a bigger compression body. 
\end{proof}

\begin{dfn}
\label{d:amalgam}
Let $H$ be a GHS and $\Sigma$ the graph in $M$ defined above. The union of the Heegaard splittings that each component of $\Sigma$ is a spine of is called the {\it amalgamation} of $H$ and will be denoted $\amlg{H}$.
\end{dfn}

Note that although the construction of the graph $\Sigma$ involved some choices, its neighborhood $N(\Sigma)$ is uniquely defined.

%\begin{lem}
%\label{l:AmalgGenus}
%Suppose $M$ is irreducible, $H$ is a GHS of $M$ and $G$ is obtained from $H$ by weak reduction. If $g(H)=g(G)$ then $\amlg{H}$ is isotopic to $\amlg{G}$. If $g(H)>g(G)$ then $\amlg{H}$ is obtained from $\amlg{G}$ by stabilization.
%\end{lem}

\begin{lem}
\label{l:AmalgGenus}
Suppose $M$ is irreducible, $H$ is a GHS of $M$ and $G$ is obtained from $H$ by a weak reduction which is not a destabilization. Then $\amlg{H}$ is isotopic to $\amlg{G}$.
\end{lem}

\begin{proof}
Suppose $G$ is obtained from $H$ by the weak reduction $(D,E)$ for the thick level $H_+$, where $D$ is above $H_+$ and $E$ is below. Let $\Sigma$ be the graph associated with $H$ as defined above. In the first stage of weak reduction we replace $H_+$ in $\thick{H}$ with the components of $H_+/D$ and $H_+/E$. We also add the components of $H_+/DE$ to $\thin{H}$. Since $(D,E)$ is not a destabilization and $M$ is irreducible it follows that $H_+/DE$ did not contain any sphere components. Hence, to complete the formation of the weak reduction it remains only to remove parallel thick and thin levels. Now observe that if we postpone this step and we first form the associated graph $\Sigma '$ as above, then we end up with the same graph as we would have if we removed parallel thick and thin levels first. 

We now claim that a neighborhood $N(\Sigma')$ is isotopic to $N(\Sigma)$. Outside $M(H_+)$ we may assume $\Sigma$ and $\Sigma'$ coincide, so we focus our attention inside $M(H_+)$.

        \begin{figure}[htbp]
        \psfrag{D}{$D$}
        \psfrag{E}{$E$}
%        \psfrag{H}{$H_+$}
        \psfrag{A}{$\Wup(H_+)$}
        \psfrag{B}{$\Wdown(H_+)$}
        \psfrag{W}{$\Wdown(H_+/E)$}
        \psfrag{w}{$\Wup(H_+/D)$}
        \vspace{0 in}
        \begin{center}
       \includegraphics[width=5 in]{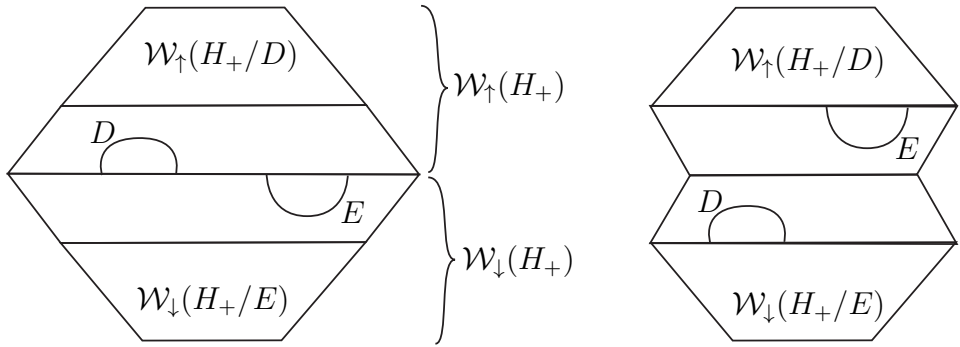}
       \caption{Because $D$ and $E$ are disjoint the manifold $M(H_+)$ can be built in two different ways.}
        \label{f:M(H_+)}
        \end{center}
        \end{figure}

The manifold $M(H_+)$ can be built as follows. Begin with $\Wdown(H_+/E)$. Attach a 1-handle to the positive boundary whose co-core is $E$. This is now the manifold $\Wdown(H_+)$. To continue, attach a 2-handle whose core is $D$, and then attach the manifold $\Wup(H_+/D)$. The 2-handle and $\Wup(H_+/D)$ is precisely $\Wup(H_+)$. See Figure \ref{f:M(H_+)}. Since $D$ and $E$ are disjoint we may build $M(H_+)$ in an alternate way: Begin with $\Wdown(H_+/E)$, attach the 2-handle, then attach the 1-handle, and finally attach $\Wup(H_+/D)$.

        \begin{figure}[htbp]
        \psfrag{D}{$D$}
        \psfrag{E}{$E$}
        \psfrag{a}{$\alpha$}
        \psfrag{g}{$\gamma$}
        \vspace{0 in}
        \begin{center}
       \includegraphics[width=5 in]{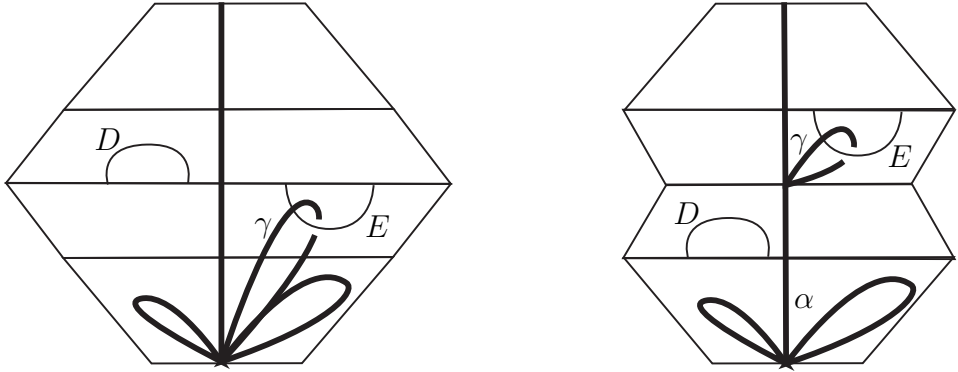}
       \caption{The graph $\Sigma$ can be obtained from $\Sigma '$ by contracting $\alpha$.}
        \label{f:Sigma}
        \end{center}
        \end{figure}

Now note that we can build $\Sigma$ by starting with a set of arcs in $\Wdown(H_+/E)$ whose complement is a product, attaching a core $\gamma$ of the 1-handle whose co-core is $E$, and then attaching some vertical arcs. Similarly, $\Sigma'$ can be built by starting with the same set of arcs in $\Wdown(H_+/E)$, attaching vertical arcs $\alpha$, then attaching $\gamma$, and finishing with more vertical arcs.  See Figure \ref{f:Sigma}. By contracting the vertical arcs $\alpha$ we achieve the  desired isotopy of $N(\Sigma')$. 
\end{proof}

\section{Sequences of Generalized Heegaard Splittings}
\label{s:SOGdefinition}

\begin{dfn}
A {\it Sequence Of GHSs} (SOG), $\{(H^i,M^i)\}$ is a finite sequence such that $H^i$ is a GHS of $M^i$ and either $H^i$ or $H^{i+1}$ is obtained from the other by a weak reduction.
\end{dfn}

\noindent {\it Notation:} We will always use superscripts to denote the GHSs of a SOG, and a boldface font to denote the entire SOG. Hence, $H^j$ is the $j$th GHS of the SOG $\bf H$.

\begin{dfn}
If $\bf H$ is a SOG and $k$ is such that $H^{k-1}$ and $H^{k+1}$ are obtained from $H^k$ by a weak reduction then we say the GHS $H^k$ is {\it maximal} in $\bf H$. 
%Similarly, if $k$ is such that $H^k$ is obtained from $H^{k-1}$ and  $H^{k+1}$ by a weak reduction then we say $H^k$ is {\it minimal}. If $H^k$ is either maximal or minimal then we say that it is {\it extremal}.
\end{dfn}

It follows from Lemma \ref{l:WRimpliesDecrease} that maximal GHSs are larger than their immediate predecessor and immediate successor. 
%Minimal GHSs are smaller then their predecessor and successor. 

\begin{dfn}
\label{d:distance}
Suppose $H$ is a Heegaard surface in a 3-manifold. Let $(V_i,W_i)$ be a reducing pair for $H$ for $i=0,1$. Then we define the {\it distance} between $(V_0,W_0)$ and $(V_1,W_1)$ to be the smallest $n$ such that there is a sequence $\{D_j\}_{j=0}^{n+1}$ where
\begin{enumerate}
	\item $\{D_0,D_1\}=\{V_0,W_0\}$,
	\item $\{D_n,D_{n+1}\}=\{V_1,W_1\}$,
	\item for all $j$ the pair $(D_j,D_{j+1})$ is a reducing pair for $H$,
	\item for $1\le j \le n$, $D_{j-1}$ is disjoint from $D_{j+1}$. 
\end{enumerate}
If there is no such sequence then we define the distance to be $\infty$.  
\end{dfn}

\begin{figure}[htbp]
        \psfrag{A}{(a)}
        \psfrag{B}{(b)}
        \psfrag{1}{$V_0$}
        \psfrag{2}{$W_0$}
        \psfrag{3}{$V_1$}
        \psfrag{4}{$W_1$}
        \psfrag{a}{$V_0=D_0$}
        \psfrag{b}{$D_2$}
        \psfrag{c}{$D_4$}
        \psfrag{d}{$D_6$}
        \psfrag{e}{$D_{1,3,5,7}$}
        \psfrag{f}{$D_9$}
        \psfrag{g}{$D_{11}$}
        \psfrag{h}{$D_{13}$}
        \psfrag{i}{$D_{15}$}
        \psfrag{j}{$D_{17}$}
        \psfrag{k}{$D_{8,10,12,14,16,18}$}
        \psfrag{l}{$D_{20}$}
        \psfrag{m}{$D_{22}$}
        \psfrag{n}{$D_{24}$}
        \psfrag{o}{$D_{26}$}
        \psfrag{p}{$D_{19,21,23,25,27}$}
        \psfrag{q}{$W_1=D_{28}$}
        \vspace{0 in}
        \begin{center}
       \includegraphics[width=3.5 in]{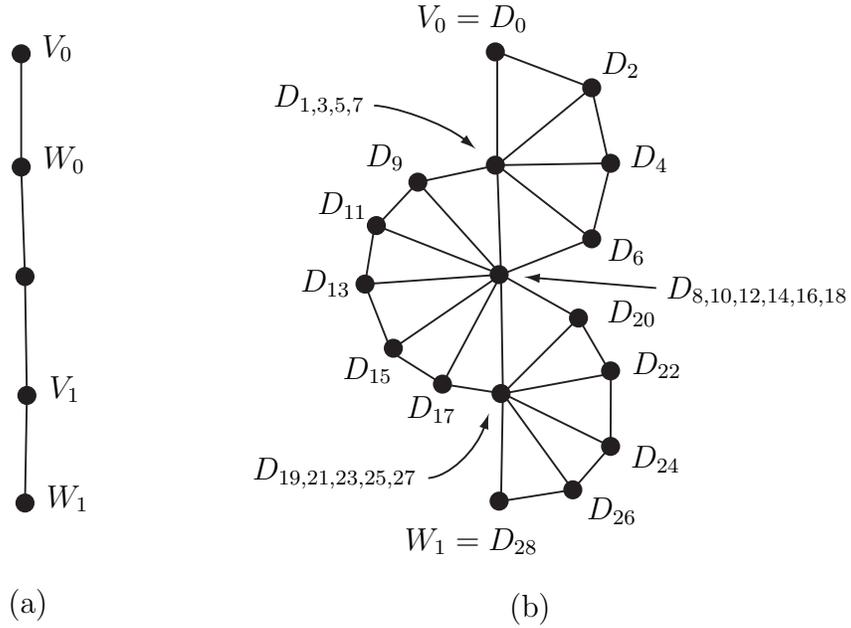}
       \caption{(a) The shortest path from $\partial V_0$ to $\partial W_1$ in the curve complex of $H$. (b) The distance from $(V_0,W_0)$ to $(V_1,W_1)$ is 27.}
       \label{f:WagonWheels}
        \end{center}
        \end{figure}

The reader may wonder how this notion of distance between reducing pairs is related to the distance you would get by dropping the last condition. This is best visualized in the curve complex of $H$, where vertices correspond to isotopy classes of essential loops in $H$ and edges correspond to disjoint pairs of loops. In Figure \ref{f:WagonWheels}(a) we have depicted a path in the curve complex from $\partial V_0$ to $\partial W_1$. In Figure \ref{f:WagonWheels}(b) we see how the boundaries of the disks $D_i$ of Definition \ref{d:distance} are related. The picture is reminiscent of the {\it geodesic hierarchies} of \cite{mm2}. The figure illustrates how the distance of Definition \ref{d:distance} is more closely related to the length of a chain of 2-simplices, rather than the length of a chain of 1-simplices. 

\begin{lem}
\label{l:FiniteImpliesFinite}
Suppose $H$ is an embedded surface in a 3-manifold. If there are reducing pairs $(V,W), (D,E)$ and $(D,E')$ for $H$ such that the distance between $(V,W)$ and $(D,E')$ is finite then the distance between $(V,W)$ and $(D,E)$ is finite.
\end{lem}

\begin{proof}
Suppose the distance between $(V,W)$ and $(D,E')$ is $n$. We now construct a sequence $\{E_i\}_{i=0}^{m}$ such that $E_0=E'$, $E_m=E$, for all $i$ the pair $(D,E_i)$ is a reducing pair, and $E_i \cap E_{i+1}=\emptyset$. The sequence 
\[\{D,E'=E_0,D,E_1,D,E_2,\dots,D,E=E_m\}\] 
then satisfies the conditions in Definition \ref{d:distance}, establishing that the distance between $(D,E')$ and $(D,E)$ is at most $2m$. The distance between $(V,W)$ and $(D,E)$ is then at most $n+2m$.

The sequence $\{E_i\}$ is defined inductively as follows:
\begin{enumerate}
	\item Define $E_0=E'$.
	\item Let $E_i$ denote the last disk defined, and suppose there is a simple closed curve in $E \cap E_i$. Let $e$ denote an innermost subdisk of $E$ bounded by a loop of $E \cap E_i$. Now surger $E_i$ along $e$. The result is a disk and a sphere. Throw away the sphere and denote the disk as $E_{i+1}$. Note that $|E \cap E_{i+1}|<|E \cap E_i|$.
	\item Let $E_i$ denote the last disk defined, and suppose there are only arcs in $E \cap E_i$. Let $e$ denote an outermost subdisk of $E$ cut off by an arc of $E \cap E_i$. 
%Note that if $|D \cap E|=1$ then $E$ is a disk and there are two choices for $e$. We may thus assume $e$ is chosen so that $e\cap D=\emptyset$. 
		Now surger $E_i$ along $e$. The result is two disks, at least one of which is a compressing disk for $H$. Call such a disk $E_{i+1}$. It follows from the fact that $e \cap D=\emptyset$ that $E_{i+1} \cap D =\emptyset$. Again, note that $|E \cap E_{i+1}|<|E \cap E_i|$.
	\item If neither of the previous two cases apply then we have arrived at a disk $E_{m-1}$ such that $(E_{m-1},D)$ is a reducing pair and $E_{m-1} \cap E=\emptyset$. Now let $E_m=E$ and we are done.
\end{enumerate}
\end{proof}

\begin{lem}
\label{l:InfiniteImpliesCritical}
Suppose $H$ is an embedded surface in a 3-manifold. If there are reducing pairs $(V_0,W_0)$ and $(V_1,W_1)$ for $H$ such that the distance between $(V_0,W_0)$ and $(V_1,W_1)$ is $\infty$ then $H$ is critical. 
\end{lem}

\begin{proof}
Let $C_0$ be the set of compressing disks such that for each $D \in C_0$ there exists an $E$ where the distance between $(V_0,W_0)$ and $(D,E)$ is finite. Let $C_1$ denote the set of compressing disks that are not in $C_0$.  We claim that the sets $C_0$ and $C_1$ satisfy the conditions of Definition \ref{d:critical}. 

Clearly, $V_0$ and $W_0$ are in $C_0$. We claim $V_1$ is in $C_1$. By symmetry it will follow that $W_1 \in C_1$, and hence Condition 1 of Definition \ref{d:critical} is satisfied. If $V_1 \notin C_1$ then there is an $E$ which forms a reducing pair with $V_1$ such that the distance between $(V_0,W_0)$ and $(V_1,E)$ is finite. But then it would follow from Lemma \ref{l:FiniteImpliesFinite} that the distance between $(V_0,W_0)$ and $(V_1,W_1)$ is finite, a contradiction. 
 
Suppose now $D \in C_0$, and $E$ is such that $(D,E)$ is a reducing pair for $H$. To establish condition 2 we must show $E \in C_0$. By definition there is an $E'$ such that the distance between $(V_0,W_0)$ and $(D,E')$ is finite. It thus follows from Lemma \ref{l:FiniteImpliesFinite} that the distance between $(V_0,W_0)$ and $(D,E)$ is finite, and hence, $E \in C_0$.  
\end{proof}

\begin{figure}[htbp]
        \psfrag{1}{$G^1$}
        \psfrag{2}{$G^2$}
        \psfrag{3}{$G^3$}
        \psfrag{4}{$G^4$}
        \psfrag{5}{$G^5$}
        \psfrag{6}{$G^6$}
        \psfrag{7}{$G^7$}
        \psfrag{A}{$(D,E)$}
        \psfrag{B}{$(D',E')$}
        \psfrag{8}{$5$}
        \vspace{0 in}
        \begin{center}
       \includegraphics[width=3.5 in]{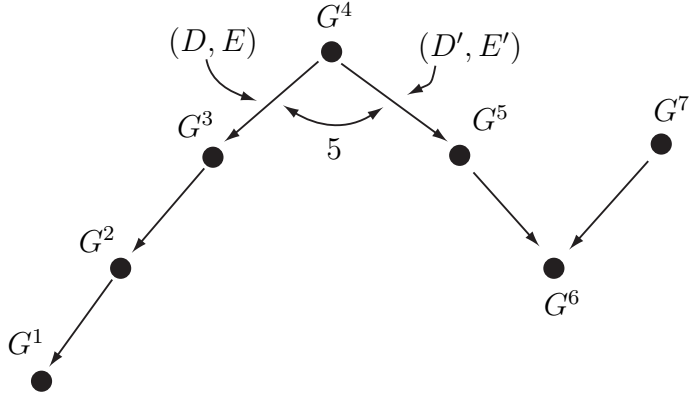}
       \caption{If the distance between $(D,E)$ and $(D',E')$ is $5$ then the angle at $G^4$ is $5$.}
       \label{f:delta}
        \end{center}
        \end{figure}

We now define a complexity on maximal GHSs of SOGs. The definition is illustrated in Figure \ref{f:delta}.

\begin{dfn}
\label{d:delta}
Suppose $G^k$ is a maximal GHS of a SOG $\bf G$. Suppose further that $G^{k-1}$ is obtained from $G^k$ by the reduction given by the reducing pair $(D,E)$ for the surface $G^k_+ \in \thick{G^k}$, and that $G^{k+1}$ is obtained from $G^k$ by the reduction given by the reducing pair $(D',E')$ for the surface $G^k_* \in \thick{G^k}$. If $G^k_+=G^k_*$ then define the {\it angle} $\angle(G^k)$ to be the distance between $(D,E)$ and $(D',E')$. Otherwise we define $\angle(G^k)$ to be 1.
\end{dfn}

Recall that a weak reduction is a way to take a GHS and obtain a smaller one. Our goal here is to define several ways to take a SOG and obtain a smaller one. Any of these will be referred to as a {\it reduction} of a SOG. To justify the statement that reduction produces something smaller, we must define a complexity for SOGs that induces a partial ordering. Furthermore, such a complexity should have the property that any decreasing sequence must terminate. This would immediately imply that any sequence of reductions must be finite. As our complexity is a lexicographically ordered multi-set of non-negative integers this latter property follows from a transfinite induction argument. 

The actual complexity we define is a bit complicated. Fortunately, the only features that we will use are easy to list:
\begin{enumerate}
	\item Eliminating or replacing a maximal GHS with one or more smaller ones will represent a decrease in complexity.
	\item Replacing a maximal GHS with angle $n$ with several identical maximal GHSs, each of which having angle less than $n$, will represent a decrease. 
%	\item Replacing a minimal GHS with a smaller one will represent a decrease (as long as this does not affect any maximal GHS).
\end{enumerate}

Any complexity one can define which behaves in this way will work for our purposes. We give one now for completeness. Let $\bf G$ be a SOG. Then the complexity of $\bf G$ is given by
\[\left\{ \left( G^k,\angle(G^k) \right) | G^k\mbox{ is a maximal GHS}\right\}\]
Sets appearing in parentheses are ordered as written. Sets appearing in brackets are put in non-increasing (lexicographical) order and repetitions are included. When comparing the complexity of two SOGs one should make lexicographical comparisons at all levels.

\begin{figure}[htbp]
        \psfrag{1}{$5$}
        \psfrag{2}{$7$}
        \psfrag{3}{$6$}
        \psfrag{4}{$5$}
        \psfrag{H}{$H$}
        \psfrag{K}{$K$}
        \psfrag{A}{}
        \psfrag{B}{}
        \psfrag{C}{}
        \vspace{0 in}
        \begin{center}
       \includegraphics[width=4.5 in]{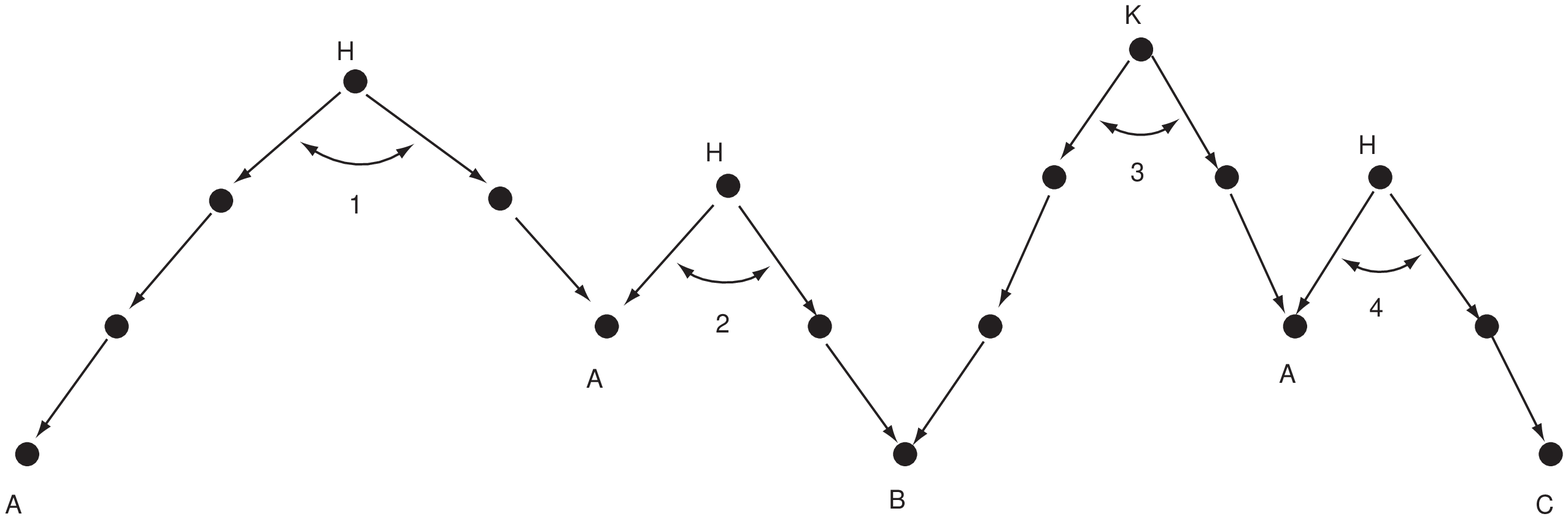}
%       \caption{If $H>K$ and $A>B>C$ then the complexity of the SOG pictured would be $(\{(H,7),(H,5),(H,5),(K,2)\},\{A,A,A,B,C\})$.}
       \caption{A SOG.}
       \label{f:MSet}
        \end{center}
        \end{figure}

\begin{ex}
Consider the SOG pictured in Figure \ref{f:MSet}. If $H>K$ then the complexity of this SOG would be \[\{(H,7),(H,5),(H,5),(K,6)\}.\] If $K>H$ then the complexity would be \[\{(K,6),(H,7),(H,5),(H,5)\}.\] Finally, if $H$ and $K$ are not comparable ({\it i.e.} $c(H)=c(K)$) then the complexity would be \[\{(H,7),(K,6),(H,5),(H,5)\}.\]
\end{ex}

We now define the various complexity decreasing operations that one can perform on a SOG that will be referred to as {\it reductions}. As in Definition \ref{d:delta} assume $G^k$ is maximal in ${\bf G}$, so that there is some thick level $G^k_+$ such that $G^{k-1}$ is obtained from $G^k$ by a reduction given by the reducing pair $(D,E)$ for the surface $G^k_+$, and there is a thick level $G^k_*$ such that $G^{k+1}$ is obtained from $G^k$ by a reduction given by the reducing pair $(D',E')$ for the surface $G^k_*$.

\subsection{Reductions of Type I}\

These are reductions which effect maximal GHSs without any consideration of their angles.

\begin{itemize}
	\item If $G^k_+ \ne G^k_*$ then we may replace $G^k$ with $H^*$ in ${\bf G}$, where $H^*$ is the GHS obtained from $G^{k-1}$ by the weak reduction $(D',E')$. Since $H^*$ can also be obtained from $G^{k+1}$ by the weak reduction $(D,E)$ our substitution has defined a smaller SOG ${\bf G}'$.  

	\item Next, assume $G^k_+=G^k_*$, but there is a thick level $G^k_0 \ne G^k_+$ in $\thick{G^k}$ which is not strongly irreducible. Let $(D^*,E^*)$ be a reducing pair for $G^k_0$. Let $H^{k-1}$, $H^k$, and $H^{k+1}$ denote the GHSs obtained from $G^{k-1}$, $G^k$, and $G^{k+1}$ by the weak reduction corresponding to $(D^*,E^*)$. Now replace $G^k$ in ${\bf G}$ with the subsequence $\{H^{k-1},H^k,H^{k+1}\}$ to define a new, smaller SOG ${\bf G}'$.
\end{itemize}

\subsection{Reductions of Type II}\

In all Type II reductions $G^k_+=G^k_*$ and we focus on $\angle(G^k)$. 

\begin{itemize}
	\item $\angle(G^k)=0$. Then $(D,E)$ is the same as $(D',E')$, so removal of the subsequence $\{G^k, G^{k+1}\}$ from ${\bf G}$ defines a new, smaller SOG ${\bf G}'$. 

	\item $\angle(G^k)=1$. In this case either $D=D'$ and $E \cap E'=\emptyset$ or $E=E'$ and $D \cap D'=\emptyset$. Assume the latter. There are now three subcases. 
	
	First, assume neither $\partial D$ nor $\partial D'$ is parallel to $\partial E$ on $G^k_+$. If $\partial D$ separates $\partial E$ from $\partial D'$ then Lemma \ref{l:WRswapSeparating} implies that $G^{k-1}$ can be obtained from $G^{k+1}$ by the weak reduction $(D,E)$. Hence, we may remove $G^k$ from $\bf G$ to obtain a smaller SOG. A symmetric argument holds if $\partial D'$ separates $\partial D$ from $\partial E$. 
	
	If neither $\partial D$ nor $\partial D'$ separates the other from $\partial E$, and neither is parallel to $\partial E$,  then we may apply Lemma \ref{l:WRswapNonSeparating}. This implies that there is a GHS $H$ that can be obtained from $G^{k-1}$ by the weak reduction $(D',E)$, and from $G^{k+1}$ by the weak reduction $(D,E)$. Hence, we may replace $G^k$ in $\bf G$ with $H$ to obtain a smaller SOG. 
	
	Finally, if $\partial D'$ is parallel to $\partial E$ we may apply Lemma \ref{l:WRswapParallel}. This implies that $G^{k+1}$ can be obtained from $G^{k-1}$ by the weak reduction $(D',E)$. Hence we may remove $G^k$ from $\bf G$ to obtain a smaller SOG. A symmetric argument holds if $\partial D$ is parallel to $\partial E$.

	\item $\angle(G^k)=n$ for some $n >1$. Let $\{D_j\}_{j=0}^{n+1}$ be a sequence given by Definition \ref{d:distance}. Choose some $m$ between 1 and $n-1$. Then the reducing pair $(D_m,D_{m+1})$ cannot be equal to either $(D,E)$ or $(D',E')$. Let $G^*$ denote the GHS obtained from $G^k$ by the weak reduction corresponding to $(D_m,D_{m+1})$. Now, let ${\bf G}'$ denote the SOG  obtained from ${\bf G}$ by inserting the subsequence $\{G^*, G^k\}$ just after $G^k$. Note that the maximal GHS $G^k$ appears one more time in ${\bf G}'$ than in ${\bf G}$. However, the angle at the old occurrence of $G^k$ is now $m$, and the angle at the new occurrence is $n-m$. As both of these numbers are smaller than $n$ we have produced a smaller SOG. 
\end{itemize}

%\subsection{Reductions of Type III}\

%A reduction of Type III is one which only affects minimal GHSs.

%\begin{itemize}
%	\item If some minimal GHS $G^k$ of $\bf G$ is not strongly irreducible then there is some weak reduction for $G^k$. Let $G^*$ denote the result of such an operation. We can now define a new SOG by inserting the subsequence $\{G^*,G^k\}$ just after $G^k$. 
%\end{itemize}

\begin{dfn}
If the first and last GHS of a SOG are strongly irreducible and none of the above reductions can be performed then the SOG is said to be {\it irreducible}.
\end{dfn}

%\begin{lem}
%\label{l:minimalGHS}
%Every minimal GHS of an irreducible SOG is strongly irreducible. 
%\end{lem}

%\begin{proof}
%The result is immediate from the fact that one cannot apply any Type III reductions.
%\end{proof}

\begin{lem}
\label{l:maximalGHS}
Every maximal GHS of an irreducible SOG is critical. 
\end{lem}

\begin{proof}
The fact that every thick level but one is strongly irreducible follows immediately from the fact that one cannot perform any Type I reductions. The remaining thick level (the surface $G^k_+$ in the definition of the reductions) must be critical by Lemma \ref{l:InfiniteImpliesCritical} since the lack of availability of Type II reductions implies that the distance between $(D,E)$ and $(D',E')$ is $\infty$.
\end{proof}

\begin{lem}
\label{l:DestabFirst}
Suppose $\bf Y$ is a SOG that is obtained from a SOG $\bf X$ by a reduction. If there is a destabilization in $\bf X$, and it comes before any stabilization, then there is a destabilization in $\bf Y$, and it comes before any stabilization in $\bf Y$. 
\end{lem}

\begin{proof}
The fact that there is a destabilization in $\bf Y$ follows immediately from the fact that no reduction will ever remove a destabilization, just possibly exchange its order with some weak reduction. The remainder of the assertion can best be seen graphically. If $X^i$ is a GHS of $\bf X$ then let $g(X^i)$ denote the genus of $\amlg{X^i}$. For the SOG $\bf X$ one can thus plot the point $(i,g(X^i))$ in the $xy$-plane. By Lemma \ref{l:AmalgGenus} the only increases and decreases in this graph are due to stabilizations and destabilizations. 

Each type of reduction effects a local maximum of this graph (although in general it will not be a strict local maximum). In each case one of the following occurs.
 	\begin{enumerate}
		\item The reduction takes a local maximum of the graph and replaces it with points at the same height or lower.
		\item The reduction removes a local maximum and possibly its successor from the graph. The successor can only be removed if it is at the same height as the predecessor. 
		\item The reduction inserts points just after a local maximum that are at the same height or lower.
	\end{enumerate} 
In each case if the first decrease (i.e. destabilization) is effected then it is replaced with a decrease. If a new increase (i.e. stabilization) is introduced then a new decrease is also introduced before it. The result follows. 
\end{proof}

%\begin{lem}
%\label{l:DestabFirst}
%Suppose $\bf G$ is a SOG such that the first destabilization comes before any stabilization (if there is one). If the SOG $\bf H$ is obtained from $\bf G$ by a reduction then the first destabilization of $\bf H$ comes before any stabilization in $\bf H$. 
%\end{lem}

%\begin{proof}
%Inspection of each type of reduction reveals the following facts:
%\begin{enumerate}
%	\item No reduction switches the order of a destabilization and a stabilization.
%	\item If a reduction introduces a stabilization then it also introduces a destabilization before it. 
%	\item If a reduction cancels a destabilization then it also cancels a stabilization that came before it. Hence, if the first destabilization comes before any stabilization then no reduction will cancel the first destabilization.
%\end{enumerate}
%The statement of the lemma follows. 
%\end{proof}

\section{An Example}

We now present an example suggested by the referee which is illustrative of some of the proof techniques used in the final section. 

Let $L_1$ and $L_2$ denote two lens spaces. Let $T_i$ denote a Heegaard torus in $L_i$. Let $\overline{T_1}$ denote the Heegaard torus obtained from $T_1$ by reversing orientation. Let $T$ denote a Heegaard torus in $S^3$.

The genus two Heegaard surfaces $T_1\# T_2$ and $\overline{T_1} \# T_2$ are equivalent after one stabilization in $M_1 \# M_2$. Hence, we may build a SOG ${\bf X}= \{(X^i,M^i)\}$ as follows:
\begin{itemize}
	\item $M^1=L_1 \cup L_2$, $\thick{X^1}=\{T_1,T_2\}$, $\thin{X^1}=\emptyset$. 
	\item $M^2=L_1 \# L_2$, $\thick{X^2}=\{T_1 \# T_2\}$, $\thin{X^2}=\emptyset$
	\item $M^3=L_1 \# L_2 \# S^3$, $\thick{X^3}=\{T_1 \# T_2 \#T\}$, $\thin{X^2}=\emptyset$
	\item $M^4=L_1 \# L_2$, $\thick{X^4}=\{\overline{T_1} \# T_2\}$, $\thin{X^2}=\emptyset$
	\item $M^5=L_1 \cup L_2$, $\thick{X^5}=\{\overline{T_1},T_2\}$, $\thin{X^1}=\emptyset$. 
\end{itemize}

The GHS $X^3$ is maximal in $\bf X$. By Lemma \ref{l:SIorCritImpliesIrreducible} it cannot be critical, since $M^3$ is reducible. Hence, by Lemma \ref{l:maximalGHS} there must be a reduction for $\bf X$. Such a reduction is not difficult to find. 

Let $S$ denote the summing sphere in $L_1 \# L_2$. Now note that $T_1$ and $\overline{T_1}$ are equivalent after one stabilization in $L_1$. Hence, when we do the connected sum with $T$ to form $X^3$ we may assume that it lies entirely on the $L_1$ side of $S$. 

The surface $T_1 \# T_2 \# T$ cuts $S$ into disks $D$ and $E$. Let $(A,B)$ denote the reducing pair which we use to go from $X^3$ to $X^2$, and let $(A',B')$ denote the reducing pair pair which we use to go from $X^3$ to $X^4$. It follows that $D$ is disjoint from both $B$ and $B'$. The fact that $B \cap B'=\emptyset$ implies the sequence $\{A,B,D,B',A'\}$ satisfies the conditions of Definition \ref{d:distance}. We conclude $\angle(X^3)=3$. (For it to be any less either $B=B'$, $A=A'$, $A \cap B'=\emptyset$, or $A' \cap B=\emptyset$. None of these are the case.) We may thus apply a reduction of Type II to $\bf G$. After a sequence of such reductions we are left with the following SOG ${\bf Y} =\{(Y^i,N^i)\}$:
\begin{itemize}
	\item $N^1=L_1 \cup L_2$, $\thick{Y^1}=\{T_1,T_2\}$, $\thin{Y^1}=\emptyset$. 
	\item $N^2=L_1 \#S^3  \cup L_2$, $\thick{Y^2}=\{T_1\# T,T_2\}$, $\thin{Y^2}=\emptyset$. 
	\item $N^3=L_1 \cup L_2$, $\thick{Y^3}=\{\overline{T _1},T_2\}$, $\thin{Y^3}=\emptyset$. 
\end{itemize}

\section{The Stability Theorem}

%In this section we combine the various results of the paper, culminating in the proof of Theorem \ref{t:mainsphere}. 

We now proceed with our proof of Gordon's conjecture.

\begin{thm}
\label{t:main}
Let $M_1$ and $M_2$ be closed, orientable 3-manifolds. Suppose $H_i$ is a Heegaard surface in $M_i$, for $i=1,2$. If $H_1 \# H_2$ is a stabilized Heegaard surface in $M_1 \# M_2$ then either $H_1$ or $H_2$ is stabilized. 
\end{thm}

\begin{proof}
Let $M_1$ and $M_2$ be two closed, orientable 3-manifolds. Let $M=M_1 \# M_2$. Suppose $H_i$ is an unstabilized Heegaard surface in $M_i$ and let $H=H_1 \# H_2$. By way of contradiction we assume there is a Heegaard surface $G$ in $M$ such that $H$ is a stabilization of $G$. Let $X$ denote the GHS of $M_1 \cup M_2$ such that $\thick{X}=\{H_1,H_2\}$ and $\thin{X}=\emptyset$. 

Let $\{M_i^j\}$ denote the irreducible manifolds in a prime decomposition of $M_i$. (If $M_i$ is the connected sum of copies of $S^2 \times S^1$ then $\{M_i^j\}=\emptyset$.) By \cite{haken:68} $H_i$ is the connected sum of Heegaard splittings $H_i^j$ of $M_i^j$ and Heegaard splittings of copies of $S^2 \times S^1$. It follows that there is a GHS $X^1$ of $\bigcup M^j_i$ such that $\thick{X^1}=\{H^j_i\}$, $\thin{X^1}=\emptyset$, and $X^1$ is obtained from $X$ by weak reduction. (Heegaard splittings of $S^2 \times S^1$ that crop up during weak reduction quickly disappear. See Example \ref{e:S^2xS^1}.) Similarly, the Heegaard surface $G$ is a connected sum of Heegaard splittings $G_i^j$ of $M^j_i$ and Heegaard splittings of copies of $S^2 \times S^1$. 

The first step is to build a SOG ${\bf X}=\{(X^i, M^i)\}_{i=1}^n$ as follows:

\begin{itemize}
	\item  $M^1=\bigcup M^j_i$,  $\thick{X^1}=\{H_i^j\}$, $\thin{X^1}=\emptyset$.
	\item $M^{i_1}=M_1 \cup M_2$, $\thick{X^{i_1}}=\{H_1,H_2\}$, $\thin{X^{i_1}}=\emptyset$.
	\item $M^{i_2}=M$, $\thick{X^{i_2}}=\{H\}$, $\thin{X^{i_2}}=\emptyset$.
	\item $M^{i_3}=M$, $\thick{X^{i_3}}=\{G\}$, $\thin{X^{i_3}}=\emptyset$. 
	\item $M^n=\bigcup M^j_i$ , $\thick{X^n}=\{G_i^j\}$, $\thin{X^n}=\emptyset$. 
\end{itemize}
For $i <i_2$ in the above SOG the GHS $X^i$ is obtained from $X^{i+1}$ by weak reduction. For $i \ge i_2$ the GHS $X^{i+1}$ is obtained from $X^i$ by weak reduction. Hence $X^{i_2}$ is maximal. 

Apply a maximal sequence of weak reductions to $X^1$ and $X^n$ to obtain strongly irreducible GHSs $X^-$ and $X^+$. We now extend $\bf X$ in the natural way to a SOG $\overline{\bf X}$ such that:
\begin{itemize}
	\item The first GHS of $\overline {\bf X}$ is $X^-$ and the last GHS is $X^+$. 
	\item The SOG $\overline {\bf X}$ contains $\bf X$ as a subsequence.
\end{itemize}

Now apply a maximal sequence of reductions to $\overline{\bf X}$ to obtain an irreducible SOG ${\bf Y}=\{Y^i\}_{i=1}^m$. By Lemma \ref{l:maximalGHS} the maximal GHSs of $\bf Y$ are critical. By Lemma \ref{l:SIorCritImpliesIrreducible} the manifolds that these are GHSs of are irreducible. It follows that every GHS in $\bf Y$ is a GHS of an irreducible manifold. 

Since $X^{i_1}$ is unstabilized (by assumption) it follows that $X^1$ is unstabilized. Since $X^1$ contains a single thick level in each $M_i^j$ it immediately follows that  $\amlg{X^1}=X^1$, and hence $\amlg{X^1}$ is unstabilized. The GHS $X^-$ is obtained from $X^1$, a GHS of an irreducible manifold, by weak reductions. As $X^1$ is unstabilized these weak reductions can not be destabilizations. Hence, by Lemma  \ref{l:AmalgGenus} we conclude $\amlg{X^-}=\amlg{X^1}$, and so $\amlg{X^-}$ is unstabilized. Finally, since $X^-$ is strongly irreducible, and $\bf Y$ is obtained from $\bf X$ by reductions, $Y^1=X^-$. Hence $\amlg{Y^1}$ is unstabilized. 

By construction there is a destabilization in $\bf X$, and this comes before any stabilization. By Lemma \ref{l:DestabFirst} there is thus an $i$ such that $Y^{i+1}$ is obtained from $Y^i$ by destabilization. Also by Lemma \ref{l:DestabFirst} if $Y^{j+1}$ is obtained from $Y^j$ by stabilization then $j>i$. It now follows from Lemma \ref{l:AmalgGenus} that $\amlg{Y^{k-1}}=\amlg{Y^k}$ for all $k \le i$. In particular, $\amlg{Y^i}=\amlg{Y^1}$, and hence $\amlg{Y^i}$ is unstabilized. We have now reached a contradiction, since $Y^{i+1}$ is obtained from $Y^i$ by destabilization.
\end{proof}

\begin{thm}
\label{t:IsotopyTheorem}
Let $M_1$ and $M_2$ be two closed, orientable 3-manifolds. Suppose $H_i $ and $G_i$ are non-isotopic Heegaard splittings of $M_i$ and $H_i$ is unstabilized. Then $H_1 \# H_2$ is not isotopic to $G_1 \# G_2$ in $M_1 \# M_2$. 
\end{thm}

\begin{proof}
The proof is similar to that of Theorem \ref{t:main}. The main difference is that the Heegaard splittings $G_1$ and $G_2$ are now given in the hypotheses of the theorem. Also, since there will be no destabilizations in the SOG $\overline {\bf X}$ the final contradiction is slightly different. 

Let $\{M_i^j\}$ denote the irreducible manifolds in a prime decomposition of $M_i$. By \cite{haken:68} $H_i$ is the connected sum of Heegaard splittings $H_i^j$ of $M_i^j$ and Heegaard splittings of copies of $S^2 \times S^1$. Also by \cite{haken:68} the Heegaard surface $G_i$ is a connected sum of Heegaard splittings $G_i^j$ of $M^j_i$ and Heegaard splittings of copies of $S^2 \times S^1$. 

By way of contradiction, we now assume $H_1 \# H_2$ is isotopic to  $G_1 \# G_2$ and build a SOG $\bf X$ as follows:

\begin{itemize}
	\item  $M^1=\bigcup M^j_i$,  $\thick{X^1}=\{H_i^j\}$, $\thin{X^1}=\emptyset$.
	\item $M^{i_1}=M_1 \cup M_2$, $\thick{X^{i_1}}=\{H_1,H_2\}$, $\thin{X^{i_1}}=\emptyset$.
	\item $M^{i_2}=M$, $\thick{X^{i_2}}=\{H_1 \# H_2 \simeq G_1 \# G_2\}$, $\thin{X^{i_2}}=\emptyset$.
	\item $M^{i_3}=M_1 \cup M_2$, $\thick{X^{i_3}}=\{G_1,G_2\}$, $\thin{X^{i_3}}=\emptyset$.
	\item $M^n=\bigcup M^j_i$ , $\thick{X^n}=\{G_i^j\}$, $\thin{X^n}=\emptyset$. 
\end{itemize}
For $i <i_2$ in the above SOG the GHS $X^i$ is obtained from $X^{i+1}$ by weak reduction. For $i \ge i_2$ the GHS $X^{i+1}$ is obtained from $X^i$ by weak reduction. 

Apply a maximal sequence of weak reductions to $X^1$ and $X^n$ to obtain strongly irreducible GHSs $X^-$ and $X^+$. We now extend $\bf X$ in the natural way to a SOG $\overline{\bf X}$ such that:
\begin{itemize}
	\item The SOG $\overline {\bf X}$ contains $\bf X$ as a subsequence.
	\item The first GHS of $\overline {\bf X}$ is $X^-$ and the last GHS is $X^+$. 
\end{itemize}

Now apply a maximal sequence of reductions to $\overline{\bf X}$ to obtain an irreducible SOG ${\bf Y}=\{Y^i\}_{i=1}^m$. Note that $Y^1=X^-$ and $Y^m=X^+$. As in the proof of Theorem \ref{t:main} it follows from Lemmas \ref{l:maximalGHS} and \ref{l:SIorCritImpliesIrreducible} that every GHS in $\bf Y$ is a GHS of an irreducible manifold. 

Since $X^{i_1}$ is unstabilized (by assumption) it follows that $X^1$ is unstabilized. Since $X^1$ contains a single thick level in each $M_i^j$ it immediately follows that  $\amlg{X^1}=X^1$. The GHS $X^-$ is obtained from $X^1$, a GHS of an irreducible manifold, by weak reductions that are not destabilizations. By Lemma  \ref{l:AmalgGenus} we conclude $\amlg{X^-}=\amlg{X^1}$. Finally, since $X^-$ is strongly irreducible, and $\bf Y$ is obtained from $\bf X$ by reductions, $Y^1=X^-$. Hence $\amlg{Y^1}=\amlg{X^-}=\amlg{X_1}=X^1$.

By construction there are no stabilizations in $\bf X$. Hence, by Lemma \ref{l:DestabFirst} if there is a stabilization in $\bf Y$ it comes after a destabilization. Suppose now there is an $i$ such that $Y^{i+1}$ is obtained from $Y^i$ by destabilization. 
%By Lemma \ref{l:DestabFirst},  if $Y^{j+1}$ is obtained from $Y^j$ by stabilization then $j>i$. 
It now follows from Lemma \ref{l:AmalgGenus} that $\amlg{Y^{k-1}}=\amlg{Y^k}$ for all $k \le i$. In particular, $\amlg{Y^i}=\amlg{Y^1}$, and hence $\amlg{Y^i}$ is unstabilized. We have now reached a contradiction, since $Y^{i+1}$ is obtained from $Y^i$ by destabilization. We conclude that there are no destabilizations, and hence no stabilizations, in $\bf Y$. But now it follows from Lemma \ref{l:DestabFirst} that there were no destabilizations in $\bf X$. We may thus conclude, as above, that $\amlg{Y^m}=\amlg{X^+}=\amlg{X^{n}}=X^{n}$.

It also follows from Lemma \ref{l:AmalgGenus} that $\amlg{Y^i}$ is the same for all $i$, and hence $\amlg{Y^1}=\amlg{Y^m}$. Finally, this gives us $X^1=X^{n}$, and the result follows. 
\end{proof}

\bibliographystyle{alpha}
\bibliography{StabilityII}

\end{document}